\documentclass{article}
\usepackage[a4paper, total={6in, 8in}]{geometry}
\usepackage[utf8]{inputenc}
\usepackage[T1]{fontenc}
\usepackage{color}
\usepackage{amsmath,bm}
\usepackage{amsthm}
\usepackage{amssymb}
\usepackage{graphicx}
\usepackage{subcaption}
\usepackage{multirow}
\usepackage{xcolor}
\usepackage{comment}
\usepackage{hyperref}
\usepackage{enumerate}
\makeatletter
\usepackage{pifont}
\newcommand{\xmark}{\ding{55}}
\newtheorem{assumption}{Assumption}
\theoremstyle{plain}

\theoremstyle{plain}
\newtheorem{thm}{\protect\theoremname}
\newtheorem{remark}{Remark}
\usepackage{algorithm}
\usepackage{cite}

\usepackage{amsmath}
\usepackage{algorithmic}

\makeatother

\newtheorem{lem}[thm]{\protect\lemmaname}

\providecommand{\lemmaname}{Lemma}
\providecommand{\propositionname}{Proposition}
\providecommand{\remarkname}{Remark}
\providecommand{\theoremname}{Theorem}

\providecommand{\propositionname}{Proposition}
\providecommand{\remarkname}{Remark}
\providecommand{\theoremname}{Theorem}

\newcommand{\vb}{\mathbf{v}}
\newcommand{\xb}{\mathbf{x}}
\newcommand{\yb}{\mathbf{y}}
\newcommand{\Q}{\textbf{Q}'}

\newcommand{\ka}{\tau}
\newcommand{\ta}{\tau}
\newcommand{\la}{\bar{\lambda}}
\newcommand{\1}{1}
\newcommand{\n}{\nabla}
\newcommand{\x}{\bar{x}}
\newcommand{\EX}{\mathbb{E}}

\newcommand{\dd}{\Delta x^*}

\newcommand{\ep}{\epsilon}
\newcommand{\xx}{\mathbf{\bar{x}}}
\newcommand{\N}{\mathbb{N}}
\newcommand{\defeq}{\stackrel{\rm def}{=}}
\newcommand{\II}{\textbf{I}_n}
\newcommand{\III}{\bar{\textbf{I}}_n}
\newcommand{\Eb}{\textbf{E}}
\newcommand{\g}{\textbf{g}}
\newcommand{\ddd}{\textbf{d}}
\newcommand{\ty}{\texttt}
\newcommand{\rb}[1]{{\color{black}#1}}

\usepackage{authblk}
\usepackage{hyperref}
\title{A Stochastic Gradient Tracking Algorithm  for Decentralized Optimization With Inexact Communication}
\author{Suhail M. Shah$^*$  \qquad \text{  }\qquad\qquad   Raghu Bollapragada\thanks{University of Texas at Austin. \texttt{(\href{mailto:suhail.mohmad@utexas.edu}{suhail.mohmad@utexas.edu}, \href{mailto:raghu.bollapragada@utexas.edu}{raghu.bollapragada@utexas.edu})}}}

\begin{document}
\maketitle
\begin{abstract}
Decentralized optimization is typically studied under the assumption of noise-free transmission. However, real-world scenarios often involve the presence of noise due to factors such as additive white Gaussian noise channels or probabilistic quantization of transmitted data.
These sources of noise have the potential to degrade the performance of decentralized optimization algorithms if not effectively addressed. In this paper, we focus on the noisy communication setting and propose an algorithm that bridges the performance gap caused by communication noise while also mitigating other challenges like data heterogeneity.
We establish theoretical results of the proposed algorithm that quantify the effect of communication noise and gradient noise on the performance of the algorithm. 
Notably, our algorithm achieves the optimal convergence rate for minimizing strongly convex, smooth functions in the context of inexact communication and stochastic gradients. Finally, we illustrate the superior performance of the proposed algorithm compared to its state-of-the-art counterparts on machine learning problems using MNIST and CIFAR-10 datasets.

\end{abstract}
\section{Introduction}
The seminal works \cite{tsit1,tsit2}  were one of the earliest works to formally study the problem of decentralized decision making and optimization. These works helped launch the field of decentralized optimization, where a connected network of multi agents collectively optimize an objective function by only exchanging information between neighboring agents in the network. 
Over the last four decades, this area has intermittently experienced phases of extensive research activity with the current iteration being mainly spurred by machine learning (ML) based optimization on decentralized data among other applications. Adding a distributed component to an optimization algorithm for ML naturally lends itself to several advantages over its centralized counterparts such as data privacy and fault tolerance while improving scalability with problem size. Formally, the problem of decentralized optimization in its most succinct form can be stated as:
\begin{align}\label{mainprob}
\min_{x_i \in \mathbb{R}^d} \,&\textbf{f}(\textbf{x}) \defeq \frac{1}{n} \sum_{i=1}^n f_i(x_i)\nonumber\\
\text{s.t.} \, & x_i =x_j,\,\forall\, i,j \in \{1,2,\cdots,n\}
\end{align}
where $\textbf{x} := (x_1, \cdots,x_n) \in \mathbb{R}^{nd} $ with $x_i$ being the copy of the optimization variable held by the $i$th node (agent) of a network and $f_i:\mathbb{R}^d \to \mathbb{R}$ is the expected value $ f_i(.) = \EX_{\xi_i}\, [  F_i(.,\xi_i)]$ of the stochastic function $F_i(.,\xi_i):\mathbb{R}^d \to \mathbb{R}$ 
private to node $i$. Problems of this nature arise in several applications with a prominent example being machine learning, where the $f_i$ is a function of the data held at node $i$.

A key aspect of decentralized algorithms is the need for communication between nodes to achieve consensus ($x_i =x_j,\,\forall\, i,j \in [n]:=\{1,2,\cdots,n\}$). However, this communication is typically not noise-free, and any form of inexactness in the algorithm can potentially degrade its performance if not addressed properly. Even 
fundamental algorithms like decentralized gradient descent (DGD) do not 
possess convergence guarantees or assured performance in the presence of inexact communication \cite[Theorem III.8]{near} or, Section IV, \textit{ibid}. Therefore, it is essential to develop a framework that incorporates inexact communication to design algorithms that effectively mitigate its adverse effects.

Data heterogeneity poses another challenge in decentralized optimization where the training data is decentralized over the nodes or generated on client devices so that each node has only access to $f_i(\cdot)$.  
Fundamental algorithms such as stochastic decentralized  gradient descent (\texttt{S-DGD}), used to solve (\ref{mainprob}) are adversely affected by data heterogeneity \cite{sgt1}. To overcome these limitations, Gradient Tracking (\texttt{GT}) type methods \cite{EXTRA,DIG} have been developed which communicate an additional vector that tracks the gradient of the global objective function. However, any inexactness in the communication can again severely degrade the overall performance \cite{near, yuan2}. In fact, with quantization, \texttt{GT} can empirically show divergent behaviour \cite[Section IV]{near}. 

In this paper we consider the question of whether the inadequacies in performance resulting from inexact communication in decentralized algorithms can be properly
addressed while retaining the benefits such as achieving consensus or removing data heterogeneity dependence. 
Specifically, our focus 
is on designing and analyzing algorithms based on the \texttt{GT} strategy in the setting where the information, which could be the current iterate or the gradient tracking vector, is corrupted by additive zero-mean noise with finite variance.

\begin{table}\label{tbt}
\footnotesize
\caption{\small{Comparison of convergence rates for strongly convex, smooth functions with stochastic gradients/communication noise for related works.} }
\centering
\begin{tabular}{ c  c  c c  }
\hline 
\hline
Reference     & Grad. Noise & Comm. Noise &  No. of iterations to $\epsilon$-acc.\\
\hline 
\hline
\cite{DIG,gt,sgt1} - \texttt{Gradient Tracking (GT)} &\xmark&\xmark& $\mathcal{O}\Big(\frac{L}{\mu\ka} \log \frac{1}{\epsilon}\Big)$\\
\hline
\cite{sgt2} - \texttt{Stochastic  DGD} &\checkmark&\xmark& $\mathcal{O}\Big(\frac{1}{n \mu }\frac{\sigma_g^2}{\epsilon}+\frac{\sqrt{L}}{\mu\ka} \frac{\sigma_g}{\sqrt{\epsilon}}   +\frac{\sqrt{L}}{\mu\ka } \frac{\chi^2}{\sqrt{\epsilon}}+ \frac{L}{\mu\ka} \log \frac{1}{\epsilon} \Big)$\\
  \hline
\cite{pudi,sgt1} -\texttt{ Stochastic GT}  &\checkmark&\xmark& $\mathcal{O}\Big(\frac{1}{n \mu }\frac{\sigma_g^2}{\epsilon}+\frac{\sqrt{L}}{\mu\ka} \frac{\sigma_g}{\sqrt{\epsilon}}  + \frac{L}{\mu\ka} \log \frac{1}{\epsilon} \Big)$\\
  \hline

\cite{com1} - \texttt{QDGD} &\xmark&\checkmark& $\mathcal{O}\Big(\frac{L^2}{\mu^2\ka} \frac{n\chi^2}{\epsilon^{2}}+\frac{L^2}{\mu^2\ka} \frac{n\sigma_c^2}{\epsilon^{2}}\Big)$\\
  \hline
\cite{do4} - \texttt{S-Near DGD}$^t$  &\checkmark&\checkmark & Non convergent. \\

  \hline
  \hline
\texttt{This work,} (\texttt{IC-GT}) &\checkmark&\checkmark& $\mathcal{O}\Big(\frac{1}{n \mu }\frac{\sigma_g^2}{\epsilon}+\frac{\sqrt{L}}{\mu\ka} \frac{\sigma_g}{\sqrt{\epsilon}}  +  \frac{1}{\mu\ka}  \frac{\sigma^2_c}{   \epsilon} +\frac{L}{\mu\ka} \log \frac{1}{\epsilon}\Big)$\\
  \hline
  \hline
\end{tabular}
    \caption*{\\\small{ Notation: $\sigma^2_g$: Gradient noise Variance, $\sigma^2_c$: Communication noise variance, $\chi^2$: Data heterogeneity constant satisfying $n^{-1} \sum_{i=1}^n \|\n f_i(x^*)- \n f(x^*)\|^2 \leq \chi^2 $ for optimal point $x^*$. 
    
$L,\mu,\ka,n$: Smoothness constant, strong convexity parameter, constant depending on network topology, total number of nodes. 

For S-Near \texttt{DGD}, $t$ denotes the number of consensus steps during each iteration and convergence is inexact even with $t\to \infty$.  The convergence is to  a neighbourhood of size $\mathcal{O}\big(\ka^{2t}\chi^2 + \frac{L^2}{\mu^2}  \sigma^2_c\big)$.
}}

\end{table}
\subsection{Related Work}
Several works 
explored the topic of inexact communication in the context of decentralized optimization, including \cite{doan, com1,do2_,do3, lee,doan2,do4,nsl,suh1,suh2,suh3}. Notably, one of the earliest and significant works in this setting are \cite{ned1,ned2}. The current work extends them in several ways, including the utilization of \texttt{GT} to address data heterogeneity and the assumptions about the underlying functions. These differences allow us to achieve superior theoretical and empirical convergence properties compared to contemporary works, as documented in Table 1 and discussed in Section \ref{sec:algo}.

Another related line of research to our work is that of decentralized optimization with randomized compressed communication \cite{cs1,cs2,cs3,cs4}. These works focus on iterate quantization for smooth and strongly convex deterministic optimization problems using randomized compression operators. However, there are significant distinctions 
between our work and these prior works, including differences in the underlying assumptions. Specifically, the algorithms proposed in the aforementioned works assume access to the compression error vector, which is transmitted to the receiving node for error compensation over a noiseless channel. Furthermore, the error variance is assumed to be controllable (\cite[Assumption 2]{cs1}) with the convergence performance being intricately linked to it(\cite[Theorem 1]{cs1}). In our setting, neither of these assumptions are applicable as they are violated in many practical scenarios, as discussed in Section \ref{sec:mod}.

The benefits of using  the \texttt{GT} strategy to address data heterogeneity have been extensively studied in numerous works \cite{DIG,gt,sgt1}.
In the deterministic setting, algorithms such as \texttt{EXTRA} \cite{EXTRA} achieve linear convergence for strongly convex, smooth functions. For the stochastic optimization setting (without communication noise), \cite{pudi,sgt1} demonstrate that \texttt{GT} based \texttt{DGD} is agnostic to the data heterogeneity. Furthermore, variants of \texttt{GT} such as \texttt{NEXT} \cite{NEXT} or the \texttt{D}$^2$ algorithm proposed in \cite{DD} have been shown to mitigate the effects of data heterogeneity. Other works exploring the \texttt{GT} strategy in various contexts include \cite{gt, uribe, DIG,sgt3,sgt4,sgt5,sgt6}.

\subsection{Contributions}

The main contributions can be summarized as follows:
\begin{itemize}
    \item[-] 
    We propose and analyze a novel variant of the 
    Gradient Tracking algorithm called Inexact Communication based Gradient Tracking (\texttt{IC-GT}) to address the challenges posed by communication noise and data heterogeneity. 
    Unlike previous approaches, our method not only retains the benefits of \texttt{GT} but also effectively eliminates the negative impact of inexact communication on algorithm performance through careful design interventions. 
    
    \item[-] 
    We show \texttt{IC-GT} can recover (upto logarithmic factors) the optimal convergence rate requirements of $\mathcal{O}(1/\epsilon)$ iterations required to achieve $\epsilon$-accuracy for stochastic optimization  while removing the data heterogeneity dependence even in the presence of communication noise. By extending the theory for exact communication based decentralized optimization \cite{st,sgt1}, our results improve upon the existing works which consider communication and gradient noise under similar assumptions and achieve either a worse convergence rate or inexact convergence (cf. Table 1).

    \item[-] To validate our theoretical results, we report experimental results that compare \texttt{IC-GT} with similar methods like \texttt{DGD} \cite{tsit1}, \texttt{DIGGing}\cite{DIG} and  \texttt{EXTRA}\cite{EXTRA}. Our experiments demonstrate the superior performance of \texttt{IC-GT} on logistic regression and image recognition problems on well known datasets.

\end{itemize}
The paper is organized as follows. We introduce the notation that is used through out the paper in the rest of this section. 
In Section \ref{sec:mod}, we describe the problem formulation and in Section \ref{sec:algo}, we present the proposed algorithm and its implementation. Section \ref{sec:conv} provides the convergence analysis while Section \ref{sec:num} presents the numerical evidence in its support. Future directions of research and conclusions are listed in Section 6.

\textit{Notation: } We use $\mathbb{R}$ to denote the set of real numbers and $\mathbb{N}$ to denote the set of all strictly positive integers. We use $\xb_k \in \mathbb{R}^{ nd}$ to denote the stacked version of $\{x_{i,k}\}_{i\in [n]}$, where $x_{i,k} \in \mathbb{R}^d$ is a column vector which denotes the value of the objective variable held by node $i$ at 
iteration $k$, i.e. $\xb_k:= (x_{1,k},\cdots, x_{n,k})$. We define $\bar{\textbf{{x}}}_k := \frac{1}{n}(1_n1_n^T\otimes I_d)\xb_k= \left(\frac{1}{n}\sum_{i=1}^n x_{i,k},\cdots,\frac{1}{n}\sum_{i=1}^n x_{i,k} \right) $, where the column vector $1_n:= (1,\cdots,1) \in \mathbb{R}^{n}$ and $I_d \in \mathbb{R}^{d\times d}$ being the identity matrix. The symbol $\otimes$ is used to denote the Kronecker product between any two matrices while $\|\cdot\|$ is understood to be the $\ell_2 $-norm of a vector or a matrix depending upon the argument. The $\ell_2$ inner product between any two vectors is denoted using $\langle \cdot,\cdot \rangle$. The following notation is used for the gradients,
$$
\n \textbf{f}( \xb_k) := \left(\n f_1 (x_{1,k}), \cdots, \n f_n(x_{n,k})\right) \text{ and }  \n \textbf{f}(\bar{\textbf{x}}_k) :=\left(\n f_1 (\bar{x}_k),\, \cdots\, ,\n f_n(\bar{x}_k)\right).
$$
We also define the matrices,
$$
\textbf{I}_n = I_n \otimes I_d \quad \text{ and } \quad \III := \II - \frac{1_n1_n^T \otimes I_d}{n}.
$$
Finally, for any two real valued functions $f(\cdot)$ and $g(\cdot)$, $f(x) = \mathcal{O}(g(x))$ denotes the standard Big-O notation which implies that there exists a finite constant $C>0$ and $x_0$ such that $|f(x)| \leq Cg(x)$ for all $x\geq x_0$.  We use $\mathcal{\tilde{O}}(\cdot)$ when ignoring logarithmic factors. 

\section{Preliminaries}\label{sec:mod}
In this section, we provide preliminaries regarding the network and communication model, and also state the assumptions that are used in the paper.

The network is represented by a (undirected) graph $\mathcal{G} =\{\mathcal{V},\mathcal{E}\}$, where $\mathcal{V}$ denotes the set of nodes and $\mathcal{E}$ represents the set of edges. 
We use the matrix $Q = [q_{ij}]_{i,j\in[n]} \in \mathbb{R}^{n \times n}$ to denote the mixing matrix (or consensus matrix) that captures the connectivity of the network. By this, we mean that the entry $q_{ij}>0$ (assumed to be equal to $q_{ji}$), if there is 
an edge between any two nodes $i,\,j\in \mathcal{V}$. We use $\mathcal{N}(i)$ to denote the set of neighbours of $i$, i.e., the set $j\in \mathcal{V}$ with $j\neq i$ for which $q_{ij}>0$. 
We make the following assumption regarding the matrix $Q$. 
\begin{assumption}[Mixing matrix]\label{asmp1}  
The mixing matrix $Q$ is symmetric and doubly stochastic. Furthermore, the eigenvalues $\{\lambda_i\}_{i\in [n]}$ of $Q$ satisfy $1 = \lambda_1 > \lambda_2 \geq \cdots \geq \lambda_n > -1$.
\end{assumption}
\begin{remark}
The symmetric and double stochasticity assumption of $Q$ is standard in decentralized optimization along with $\lambda_2<1$ which implies that the graph is connected. Therefore, it implies that $(Q \otimes I_d )\textbf{x } = \textbf{x }$ if and only if $x_i = x_j$ for all $i,j \in \mathcal{V}$. 
Moreover, it also ensures that the spectral gap $\delta(Q):=1- \max\{|\lambda_2|,|\lambda_n|\}$ is greater than zero which in turn ensures that the consensus error decreases linearly after each averaging step, i.e.,
\begin{equation}\label{contract}
\left\| (Q \otimes I_d ) \textbf{x }- \left(\frac{1_n1_n^T}{n}\otimes I_d \right)  \textbf{x}\right\|^2 \leq (1-\delta)^2 \left\| \textbf{x}  - \left(\frac{1_n1_n^T}{n}\otimes I_d \right)  \textbf{x}\right\|^2
\end{equation}
for any $\textbf{x}\in \mathbb{R}^{nd}$. For undirected graphs, this assumption can be guaranteed by using the Metropolis weights(\cite[Section 3]{DIG}).
\end{remark}
We next describe the communication model considered in this work. 
We make the assumption that when any node $i \in [n]$ sends a signal vector $x_{i,k} \in \mathbb{R}^d$ to a neighboring node $j$ at iteration $k \in \mathbb{N}$, node $j$ receives the vector $\varphi_c(x_{i,k}) \in \mathbb{R}^d$ instead of the original vector $x_{i,k}$ , where $\varphi_c(\cdot):\mathbb{R}^d \to \mathbb{R}^d$ represents a random transformation given by
$$\varphi_c(x_{i,k}) := x_{i,k} + \epsilon_{i,k,c}\,,$$
where $\epsilon_{i,k,c} \in \mathbb{R}^d$ is a random vector. 
We emphasize that we do not assume access to the values of $\epsilon_{i,k,c}$. We make the following assumption concerning $\epsilon_{i,k,c}$.
\begin{assumption}[Noisy signal transmission]  \label{asmp2}
The random noisy vector $\epsilon_{i,k,c}$ is assumed to be zero mean conditioned on $x_k$ with bounded variance for all $i \in [n]$ and $k\in \mathbb{N}$, i.e., 
$$
\EX \,[\epsilon_{i,k,c}|x_{i,k}]=0 ,\,\,\,\,\,\EX\,[\|\epsilon_{i,k,c}\|^2]\leq \sigma^2_{c},
$$
for some finite $\sigma_c>0$.
\end{assumption}
\noindent We describe two important examples of $\varphi_c(\cdot)$ for which Assumption \ref{asmp2} is satisfied.\\

\noindent \textbf{Additive White Gaussian Noise channel (AWGN)}: The most common approach to modeling an analog based communication channel between two nodes is through an AWGN channel \cite{AWGN}. In this scenario, when a node transmits a signal $y_{\text{tr}} \in \mathbb{R}$ to a neighboring node, the received signal at the receiving node, denoted as $y_{\text{rc}}$, can be represented as
$$
y_{\text{rc}} = h y_{\text{tr}} + \epsilon_{c},
$$
where $h \in \mathbb{R}$ captures channel effects like fading \cite{fading}, and $\epsilon_{c}$ represents zero-mean Gaussian noise with variance $\sigma_c^2$, independent of the transmitted signal $y_{\text{tr}}$. Assuming that the receiving node possesses a prior estimate of $h$ \cite[Chapter 4]{book}, it can construct an estimate of the true signal $\hat{y}_{\text{rc}}$ as, 
$$
\hat{y}_{\text{rc}} =  \frac{1}{h} y_{\text{rc}}= y_{\text{tr}}  + \frac{\epsilon_{c}}{h}.
$$
Hence, in this scenario, we can express $\varphi_c(\cdot)$ as,
$$
\varphi_c(y_{\text{tr}}) = y_{tr} + \frac{\epsilon_{c}}{h},
$$
implying Assumption \ref{asmp2} is satisfied since $\EX \,[\varphi_c(y_{\text{tr}})] = y_{\text{tr}} $ and $\EX \,[\|\varphi_c(y_{\text{tr}})- y_{\text{tr}}\|^2] \leq \sigma_c^2/h^2 $.\\

\noindent \textbf{Probabilistic Quantization:} Another significant example of operator $\varphi_c$ arises in the context of quantization with unbiased compression operators. Specifically, consider a scalar $x \in \mathbb{R}$. The quantized value $\varphi_c(x)$ can be determined based on the following rule:
\begin{equation}\label{randquant}
\varphi_c(x)= 
    \begin{cases}
    \lfloor x \rfloor_p \text{ with probability } (  \lceil x \rceil_p -x)\Delta_p\\
    \lceil x \rceil_p \text{ with probability } (x-  \lfloor x \rfloor_p )\Delta_p
    \end{cases}
\end{equation}
where $\lfloor x \rfloor_p$ and $\lceil x \rceil_p$ denote the operations of rounding down and up to the nearest integer multiple of $\frac{1}{\Delta_p}$ respectively, and $\Delta_p$ is a positive integer. The operator $\varphi_c$ defined in (\ref{randquant}) satisfies $\EX \,[\varphi_c(x)] =x$ and $\EX\,[| \varphi_c(x) -x|^2] \leq \tfrac{1}{4\Delta_p^2}$ as shown in \cite{pq} implying Assumption \ref{asmp2} is satisfied.\\

\noindent We also make the following assumptions regarding the objective function. 
 
\begin{assumption}[Regularity and convexity]\label{asmp3}  
Each local function $f_i$ is $L$-smooth and 
$\mu$-strongly convex.
\end{assumption}
\begin{assumption}[Unbiased Gradient Samples]\label{asmp4} 
Each node $i$ has access to conditionally unbiased, finite variance gradient samples $\n F_i(x_{i,k},\xi_k)$ of $\n f_i(x_{i,k})$ for any given $x_{i,k}\in \mathbb{R}^d,$ $k\in \mathbb{N}$. That is,  
$$
\EX_{\xi_{i,k}}\, [\n  F_i(x_{i,k},\xi_{i,k})\,|\,x_{i,k}] = \n f_i(x_{i,k}) \text{ and }  \EX_{\xi_{i,k}} \,[\|\n F_i(x_{i,k},\xi_{i,k}) -\n f_i(x_{i,k})\|^2] \leq \sigma^2_{g}
$$
for some finite $\sigma_g>0$ with $\xi_{i,k}$ being assumed to be independent of $\epsilon_{i,k,c}$.
\end{assumption} 
\begin{remark}
The finite variance assumption in Assumption 4 can be relaxed along two possible lines with minor modifications to the convergence analysis. One relaxation would be to allow the noise to grow with the gradient norm (cf. Assumption 3b, \cite{st}). The other possibility is to replace $\sigma^2$ with $\sigma^2_*:= \frac{1}{n} \sum_{i=1}^n \|\n F(x^*,\xi_i)- \n f(x^*)\|^2$, the noise at the optimal point $x^*$, as in \cite{ngu}. 
\end{remark}

\begin{remark}
The convergence analysis can also be extended to a non-convex setting by modifying the measure of stationary to be the $\ell_2$-norm of the gradient. 
\end{remark}

\section{The \texttt{IC-GT} method}\label{sec:algo}
In this section, we describe the proposed method that accounts for inexact communication, referred to as Inexact Communication based Gradient Tracking \ty{(IC-GT)} designed to solve the problem (\ref{mainprob}). Algorithm 1 presents the pseudo code of \ty{(IC-GT)}.
\begin{algorithm}[H]
\label{alg:code1}
\footnotesize
   \caption{\small \texttt{INEXACT COMMUNICATION based GRADIENT TRACKING (IC-GT)}}
    \footnotesize
  \begin{algorithmic}[1]
   \footnotesize
     \STATE \textbf{Input} Graph $\mathcal{G}(\mathcal{V},\mathcal{E})$; Matrix $Q= [q_{ij}]_{i,j\in [n]} \in \mathbb{R}^{n \times n}$ ; Operator $\varphi_c(\cdot)$;  Noise attenuation parameter $\gamma > 0$; Step size parameter $\alpha > 0$.
     \STATE \textbf{Initialization} $x_{i,0} \in \mathbb{R}^d, \,\forall i$; $y_{i,0} := \n F_i(x_{i,0} ,\xi_{i,0} ), \,\forall i$. 
    \WHILE{$ k \geq 1$ in parallel:}
     \FOR{ all $i\in [n],$} 
         \STATE 
          $
          v_{i,k} = (1-\gamma)x_{i,k} + \gamma q_{ii} x_{i,k} +\gamma \sum_{j\in \mathcal{N}(i)}q_{ij}\varphi_c(x_{j,k})\
          $
           \STATE 
          $
           x_{i,k+1} = v_{i,k} - \alpha y_{i,k}
          $
           \STATE 
          $
          y_{i,k+1}=  (1-\gamma)y_{i,k} + \gamma q_{ii} y_{i,k}  +\gamma \sum_{j\in \mathcal{N}(i)}q_{ij}\varphi_c(y_{j,k}) + \nabla F_i(x_{i,k+1},\xi_{i,k+1}) - \nabla F_i(x_{i,k}, \xi_{i,k} ) 
          $
       \ENDFOR  
      \STATE $k  \rightarrow k + 1$ 
   
    \ENDWHILE
  \end{algorithmic}
\end{algorithm}
 To express \texttt{IC-GT} in matrix form, we introduce the matrices $Q' := [q'_{ij}]_{i,j\in [n]}$ and $\hat{Q} :=[\hat{q}_{ij}]_{i,j\in [n]}$ defined as follows: 
 \begin{equation}\label{q'hat}
 \Q\defeq \left( I_n-Q \right)\otimes I_d \qquad  \qquad   \hat{\textbf{Q}}\defeq\left( Q-\text{diag}(Q)\right)\otimes I_d,
 \end{equation}
 where $\text{diag}(Q)$ denotes the diagonal matrix with entries $q_{ij}$ for $i=j$ and 0 otherwise. Using the communication model, $\varphi_c(x_{j,k}) = x_{j,k}  + \epsilon_{j,k,c},$ we can express the iteration for $v_{i,k}$ as follows:  
 \begin{align*}
      v_{i,k} &= \left(1-\gamma(1-q_{ii})\right)x_{i,k} +\gamma \sum_{j\in \mathcal{N}(i)}q_{ij}\varphi_c(x_{j,k})\\
      &= \left(x_{i,k}-\gamma(1-q_{ii}))x_{i,k} + \gamma \sum_{j\in \mathcal{N}(i)}q_{ij} x_{j,k}\right) +\gamma \sum_{j\in \mathcal{N}(i)}q_{ij} \epsilon_{j,k,c}.
 \end{align*}
 Performing a similar manipulation for the $y$ update, we can express \texttt{IC-GT} using (\ref{q'hat}) as follows:
\begin{align}
    \vb_{k}&= (\textbf{I}_n-\gamma \Q)\xb_{k} + \gamma\hat{\textbf{Q}}\bm{\epsilon}_{k,c} \label{mat0}\\
    \xb_{k+1} &= \vb_{k} - \alpha \yb_{k} \label{matI}\\
     \yb_{k+1}&=  (\textbf{I}_n-\gamma \Q)  \yb_{k}  + \nabla \textbf{F}(\xb_{k+1},\mathbf{\bm{\xi}}_{k+1} ) - \nabla \mathbf{F}(\xb_{k},\mathbf{\bm{\xi}}_{k} ) + \gamma  \hat{\textbf{Q}}\hat{\bm{\epsilon}}_{k,c}\label{matII}
\end{align}
where $\hat{\epsilon}_{i,k,c} := \varphi_c(y_{i,k}) -y_{i,k}$, $\bm{\epsilon}_{k,c}:= 
    \left(
\epsilon_{1,k,c} ,\cdots, \epsilon_{n,k,c} 
\right)$ and $\nabla \textbf{F}(\xb_{k},\mathbf{\bm{\xi}}_{k} ) := \big( \n F_1(x_{1,k},\xi_{1,k}), \cdots,$ $\n F_n(x_{n,k},\xi_{n,k}) \big)$. \\

\noindent 
We next discuss the main modification made to the standard \texttt{DGD} algorithm \cite{tsit1} utilized in \texttt{IC-GT} to better understand its communicating and computational capabilities.\\

\noindent \textbf{(i) Use of $\textbf{I}_n - \gamma\Q$:}  In the context of \texttt{IC-GT}, the weight matrix $\II-\gamma\Q$ is employed instead of the typical $\textbf{Q}$ used in \texttt{DGD} \cite{tsit1}. To illustrate its effectiveness in mitigating communication noise, let us examine the sequence  $\{\xb_k \}_{k\geq0}$ generated according to the recursion:
\begin{align}\label{mat01}
	\xb_k &= (\textbf{I}_n-\gamma \Q)\xb_{k-1} + \gamma \hat{\textbf{Q}}\bm{\epsilon}_{k-1,c}, 
\end{align}
where the noise term $\bm{\epsilon}_{k-1,c}$ satisfies Assumption \ref{asmp2}. The recursion in (\ref{mat01}) can be interpreted as a distributed averaging algorithm using the weight matrix $\textbf{I}_n - \gamma\Q$. Specifically, when $\gamma=1$ and $\bm{\epsilon}_{k-1,c}=0$, (\ref{mat01}) reduces to the standard distributed averaging algorithm \cite{xiao}. Next, we consider the expression for the averaged iterates $\bar{\xb}_k$ obtained by multiplying (\ref{mat01}) by $\frac{1}{n}\left( 1_n 1_n^T \otimes I_d\right)$:
\begin{align}\label{mat00}
	\xx_{k} = \xx_{k-1} -\gamma \frac{1}{n} \big(1_n1_n^T\otimes I_d\big)\hat{\textbf{Q}} \bm{\epsilon}_{k-1,c},
\end{align}
where we used $\left(1_n^T \otimes I_d \right)\left(\textbf{I}_n-\gamma \Q\right) = (1^T_n\otimes I_d)$ from Assumption \ref{asmp1}. 
Subtracting (\ref{mat00}) from (\ref{mat01}) and defining $\tilde{\textbf{Q}}:= \left(\II -n^{-1} 1_n1^T_n\otimes I_d \right) \hat{\textbf{Q}}$ and recalling $\III:= \II - \frac{1_n1_n^T \otimes I_d}{n}$, we get,
\begin{align*}
	\xb_k - \xx_k &= (\textbf{I}_n-\gamma \Q)(\xb_{k-1} - \xx_{k-1}) + \gamma\tilde{\textbf{Q}}\bm{\epsilon}_{k-1,c} \\
	&=   (\textbf{I}_n-\gamma \Q)(\xb_{k-1} - \xx_{k-1}) + \gamma\tilde{\textbf{Q}}\bm{\epsilon}_{k-1,c}  -\frac{1_n1_n^T \otimes I_d}{n}(\xb_{k-1} - \xx_{k-1} ) \\
	&=  (\III-\gamma \Q )(\xb_{k-1} - \xx_{k-1}) + \gamma\tilde{\textbf{Q}}\bm{\epsilon}_{k-1,c}, 
\end{align*}	
where the second equality is due to $\frac{1_n1_n^T \otimes I_d}{n}(\xb_{k-1} - \xx_{k-1} ) = 0$.  Applying norms and taking squares yields, 

\begin{align}
	\|\xb_k - \xx_k\|^2 &\leq  \|\III-\gamma \Q \|^2 \|\xb_{k-1} -\xx_{k-1}\|^2 + \gamma^2  \|\tilde{\textbf{Q}} \bm{\epsilon}_{k-1,c}\|^2\nonumber +2\gamma\left\langle (\III-\gamma \Q)(\xb_{k-1} - \bar{\textbf{x}}_{k-1}), \tilde{\textbf{Q}}\bm{\epsilon}_{k-1,c}\right\rangle. 
\end{align}
Using the conditional zero mean and finite variance  assumption for $\bm{\epsilon}_{k-1,c}$ (Assumption \ref{asmp2}), we get,
\begin{align*}
	\EX [\|\xb_k - \xx_k\|^2]  &\leq  \left( 1- \gamma (1-  \lambda_2)\right)^{2} \EX [\|\xb_{k-1} -\xx_{k-1}\|^2] + 2n\gamma^2 \sigma_{c}^2, 
\end{align*}
where we used $\|\III-\gamma \Q \| \leq 1- \gamma(1-\lambda_2)$ (cf. (\ref{2normbd})) and $\|\tilde{\textbf{Q}}\|^2\leq 2$. Applying the above inequality repeatedly through iteration $k=0$ yields, 
\begin{align}\label{u19}
	\EX [\|\xb_k - \xx_k\|^2]  &\leq  \left( 1- \gamma (1-  \lambda_2)\right)^{2k} \|\xb_{0} -\xx_{0}\|^2 + \frac{2 n \gamma \sigma_{c}^2}{1-\lambda_2}.
\end{align}
(\ref{u19}) unveils a fundamental trade-off between two crucial aspects: the rate of decay of the consensus error and the mitigation of the influence exerted by the communication noise variance. As the parameter $\gamma$ decreases, 
a smaller final consensus error can be achieved. However, this improvement comes at the expense of a slower convergence rate in reducing the consensus error. In view of this trade-off, the parameter $\gamma$ is referred to as the `\textit{noise attenuation}' parameter. \\ 
 
\noindent \textbf{(ii) Use of Gradient Tracking:} Another crucial feature of \texttt{IC-GT} is its ability to track gradients while accommodating inexact communication through gradient tracking. The inclusion of gradient tracking offers the advantage of making the algorithm agnostic to data heterogeneity. To elaborate,
 the number of iterations required to achieve $\epsilon$-accuracy using stochastic \texttt{DGD} depends on $\mathcal{O} \left( \frac{\sqrt{L}\chi^2} {\sqrt{\epsilon}}\right)$ \cite{st}, where $\chi$ is a constant satisfies the inequality
$$
\frac{1}{n} \sum_{i=1}^n \|\n f_i(x^*)- \n f(x^*)\|^2 \leq \chi^2, 
$$
with $x^*$ denoting the optimal solution of (\ref{mainprob}). In contrast, \texttt{IC-GT} eliminates the dependence on $\chi$ entirely and, moreover, recovers the linear convergence rate in scenarios where the variances of both the gradient and communication noise are zero.

\section{Convergence Analysis}\label{sec:conv}
In this section, we establish theoretical convergence guarantees for the proposed \texttt{IC-GT} algorithm. We build up to our main result through a series of technical lemmas which we state next. 
\subsection*{Preliminaries}
For the sake of brevity, 
we assume $\bm{\epsilon}_{k,c} = \hat{\bm{\epsilon}}_{k,c}$ in (\ref{mat0})-(\ref{matII}) for all $k\in \mathbb{N}$ without loss of generality. We begin by expressing the algorithm in terms of the difference between the variables and their corresponding averages, which we refer to as the \textit{consensus} error. To denote this, we adopt the notation  $\Delta \textbf{z }:=\textbf{z} - \Bar{\textbf{z}}$ for any variable $\textbf{z}\in \mathbb{R}^{nd}$, where $\Bar{\textbf{z}}$ denotes the average, i.e. $\Bar{\textbf{z}}:= \Big( \frac{1_n1_n^T}{n}\otimes I_d\Big)z$. We first establish a recursive relation for the consensus error. 
\begin{lem}\label{lem0}
\emph{[\textbf{Recursive relation for consensus errors}]}
Suppose 
$\bm{\epsilon}_{k,c} = \hat{\bm{\epsilon}}_{k,c}$ in (\ref{mat0})-(\ref{matII}) for all $k\in \mathbb{N}$. 
Then, the iterates generated by \texttt{IC-GT} satisfy the following recursive relation:
\begin{equation}\label{mm1}
\Psi_k = \textbf{J}_\gamma \Psi_{k-1} + \alpha \Eb_{k-1},
\end{equation}
where 

\begin{equation}\label{not1}
	\Psi_k\defeq  \begin{bmatrix}
		\Delta \vb_{k}\\
		\Delta \xb_{k}\\
		\alpha \Delta \yb_{k}
	\end{bmatrix},\qquad 
	\,\,\,
	\textbf{J}_\gamma\defeq
	\begin{bmatrix}
		\III-\gamma \Q &0 &-(\III-\gamma \Q )\\
		0 & \III-\gamma \Q  &-\III\\
		0 & 0&  \III-\gamma \Q 
	\end{bmatrix},
\end{equation}
and 
$$
\Eb_{k-1} \defeq  \frac{\gamma}{\alpha}
\begin{bmatrix}
	\tilde{\textbf{Q}}\bm{\epsilon}_{k,c} \\
	\tilde{\textbf{Q}} \bm{\epsilon}_{k-1,c} \\
	\alpha\tilde{\textbf{Q}} \bm{\epsilon}_{k-1,c} 
\end{bmatrix} 
+ 
\begin{bmatrix}
	0\\
	0\\
	\III  \left(\nabla \textbf{F}(\xb_{k},\bm{\xi}_{k} ) - \nabla \textbf{F}(\xb_{k-1},\bm{\xi}_{k-1} )\right)
\end{bmatrix}
$$
with $\III:= \left(I_n  - \frac{1_n1_n^T}{n}\right)\otimes I_d$, $\Q\defeq (I_n-Q)\otimes I_d$, $\hat{\textbf{Q}}\defeq(Q - \text{diag}(Q))\otimes I_d$ and $\tilde{\textbf{Q}}  \defeq \III \hat{\textbf{Q}} $.
\end{lem}
The proof of this lemma is provided in Appendix \hyperref[sec:apndxI]{I}. One of the challenges in analyzing \texttt{IC-GT} is that the matrix $\textbf{J}_{\gamma}$ defined in (\ref{not1}) is not necessarily a contractive matrix. In other words, the condition $\|\textbf{J}_{\gamma}\|<1$ is not guaranteed to hold. However, the following result demonstrates that despite this restriction, there exists a positive integer $\tau$ such that $\|\textbf{J}_{\gamma}\|^\tau <1$.
\begin{lem}\label{lem3}
\emph{[\textbf{Strict contractive property for $\textbf{J}_\gamma$}]}
Suppose Assumption~\ref{asmp1} holds. 
For any given 
$\delta \in (0,1)$, $\gamma \in (0,1/4)$ and $\lambda_2 $ associated with the matrix $Q$, suppose $\ka\in \mathbb{N}$ satisfies
\begin{equation}\label{kappaprop}
	\ka \geq \left\lceil\frac{2}{\gamma(1-\lambda_2)}\max\left\{4\ln  \left(\frac{2}{\gamma(1-\lambda_2)}\right)\,,\, \left( \gamma (1-\lambda_2) -\ln  \frac{\sqrt{\delta}}{4}  \right)\right\}\right\rceil,
\end{equation}
where $\lceil \cdot \rceil$ denotes the ceiling function. Then, $\|\textbf{J}^{\ka}_\gamma\|^2 \leq \delta<1$, where $\textbf{J}^\ka_\gamma := \underbrace{\textbf{J}_\gamma \times \cdots\times \textbf{J}_\gamma}_{\ka \text{ times }} $.
\end{lem}
The proof of this lemma is provided in Appendix \hyperref[sec:apndxII]{II}. The next result establishes a  descent relation for the consensus error $\EX \,[\|\Psi_{t+\ta}\|^2]$ in terms of $\EX [\|\Psi_{t}\|^2]$.
\begin{lem}\label{r2}
\emph{[\textbf{Descent relation for consensus error, $\EX [\|\Psi_{t}\|^2] \,$}]}
Suppose Assumptions \ref{asmp1}-\ref{asmp4} hold and $\bm{\epsilon}_{k,c} = \hat{\bm{\epsilon}}_{k,c}$ in (\ref{mat0})-(\ref{matII}) for all $k\in \mathbb{N}$. If $\gamma $ and $\alpha$ satisfy (\ref{alphacond}), then, for a given $0<\rho' \leq 1/4$, there exists a $\ka\in \mathbb{N }$ such that the following relations are satisfied for $t\geq \ta$:
\begin{align}\label{main_rec}
\EX [\|\Psi_{t}\|^2] &\leq  \rho' \EX [\|\Psi_{t-\tau}\|^2]+ 576 \alpha^2 \ka L^2  \sum_{i=t-\tau}^{t-1} \EX [\|\Psi_{i}\|^2] 
+ 1344  \alpha^2 \ka \sum_{i=t-\tau}^{t-1} \EX \left[\left\|\n \textbf{f}(\bar{\textbf{x}}_{i})-\n \textbf{f}(\textbf{x}^*) \right\|^2\right] \nonumber \\
  &\qquad + 64\gamma^2 \big(2 + \alpha^2 (\ka^2 +1/2) + \alpha^2 t\big)n\sigma_{c}^2 \ka+196n(\ka+1)\alpha^2\sigma^2_{g}
\end{align}
and for any $\ell < \ta$: 
\begin{align}\label{main_rec_2}
\EX [\|\Psi_{\ell}\|^2] &\leq 2(1+\tau^2) \|\Psi_{0}\|^2+ 576 \alpha^2 \ka L^2  \sum_{i=0}^{\ell-1} \EX [\|\Psi_{i}\|^2] 
+ 1344  \alpha^2 \ka \sum_{i=0}^{\ell-1} \EX \left[\left\|\n \textbf{f}(\bar{\textbf{x}}_{i})-\n \textbf{f}(\textbf{x}^*) \right\|^2\right] \nonumber \\
  &\qquad + 64\gamma^2 \big(2 + \alpha^2 (\ka^2 +1/2) + \alpha^2 \ell \big)n\sigma_{c}^2 \ka+196n(\ka+1)\alpha^2\sigma^2_{g}.
\end{align}
\end{lem}
The proof of this lemma is provided in Appendix \hyperref[sec:apndxII]{III}. 
We next prove an auxiliary result that will be useful for bounding the consensus error. 

\begin{lem}\label{lem4}
Suppose the non-negative scalar sequences $\{a_t\}_{t\geq0}$ and $\{e_t\}_{t\geq0}$ 
satisfy the following recursive relation for a fixed $\ka\in \mathbb{N}$:
\begin{align}
a_{t} &\leq 
\begin{cases}
    \rho' a_{t-\ka} + \frac{b}{\ka} \sum_{i=t-\ta}^{t-1} a_{i}  + c  \sum_{i=t-\ta}^{t-1} e_{i} + r  & \textit{if }t\geq \ta \label{rel} \vspace{1em} \\
    \rho'' a_{0} + \frac{b}{\ka} \sum_{i=0}^{t -1} a_{i} + c  \sum_{i=0}^{t-1}  e_{i} + r & \textit{if }t < \ta 
\end{cases}
\end{align}
where $b,\,c,\,r,\,\rho''$ are non-negative constants satisfying $b \leq \rho'/4$ and $\rho'\in \left(0,1/4\right]$. Then, for any $t \in \mathbb{N}$,
\begin{align}\label{lem5_main}
a_t &\leq 20\rho''\Big(1- \frac{3\rho}{4\ka}\Big)^t a_0 + 60 c\sum_{i=0}^{t-1}   \Big(1- \frac{3\rho}{4\ka}\Big)^{t-i}e_i +\frac{26 r}{\rho},
\end{align}
where $\rho := 1-2\rho'$.
\end{lem}
The proof of this lemma is provided in Appendix \hyperref[sec:apndxIV]{IV}. We are ready to state and prove the main convergence result. 

\subsection*{Main Result}
For convenience, we define $\dd_k$ as
\begin{align}
\dd_k \defeq \EX \left[\|\bar{x}_{k} -x^*\|^2\right], \qquad \forall k \in \N. 
\end{align}
where $x^*$ is the optimal solution of (\ref{mainprob}). 

\begin{thm}\label{thm1}
\emph{[\textbf{Convergence rate of IC-GT}]}
Suppose Assumptions \ref{asmp1}-\ref{asmp4} hold and $\bm{\epsilon}_{k,c} = \hat{\bm{\epsilon}}_{k,c}$ in (\ref{mat0})-(\ref{matII}) for all $k\in \mathbb{N}$. 
If
\begin{equation}\label{alphacond}
    \alpha  \leq \min\left\{1, \frac{1}{161280 \ka L}\right\}  \quad \mbox{and} \quad 0<\gamma <1/4,
\end{equation}
where 
\begin{equation}\label{kappabd}
\ka = \left\lceil\tfrac{2}{\gamma(1 - \lambda_2)}\max\left\{4\ln  \left(\tfrac{2}{\gamma (1-\lambda_2)}\right)\,,\, \gamma (1-\lambda_2) + \ln 16 \right\}\right\rceil.
\end{equation}
Then, for any $T\in \mathbb{N}$, we have, 
\begin{align}\label{finalresult}
\dd_{T} &\leq  (1-\alpha\mu/4 )^T\left(\dd_0 +  \frac{800(1+\ka^2)L}{n (1-\alpha\mu/4)\mu}\| \Psi_{0}\|^2 \right)  + \left( \frac{4\alpha}{\mu} + \frac{101920\,  L (\ka+1) n\alpha^2}{ \mu}\right)\frac{\sigma_g^2}{n}\nonumber\\
&+ \left( \frac{4(1+2\mu^{-1}T\alpha)}{\mu} \frac{\gamma^2}{\alpha}+ \frac{33280(2 + \alpha^2 (\ka^2+\frac{1}{2}) + \alpha^2T)L}{ \mu}n \ka \gamma^2 \right)\frac{\sigma^2_{c}}{n}.  \end{align}
\end{thm}
We make the following remarks regrading Theorem \ref{thm1}. 

\begin{remark}
  \text{(\textbf{Dependence of $\ka$ on network})} 
The parameter $\ka$ depends on the network connectivity $(\lambda_2)$ and the noise attenuation parameter $\gamma$ (cf. \ref{kappabd}) which highlights the role played by $\gamma$ in shaping the consensus properties of \texttt{IC-GT} (cf. Lemma \ref{r2}). From (\ref{kappabd}), we note that a smaller value of $\gamma$ increases $\tau$ but reduces the impact of the communication noise variance $\sigma^2_c$ in (\ref{finalresult}) which is reminiscent of the trade-off discussed in Section \ref{sec:algo}.
\end{remark}

\begin{remark} \textbf{(Iteration complexity of \texttt{IC-GT})} (\ref{alphacond}) and (\ref{kappabd}) suggest that the choices of the step size $\alpha$ and the noise attenuation parameter $\gamma$ are inherently connected. Using (\ref{kappabd}) in (\ref{alphacond}), we have the following relation:
\begin{equation}\label{50}
\frac{\alpha}{\gamma} = 
\tilde{\mathcal{O}} \left(\frac{1-\lambda_2}{L}\right)
\end{equation}
To calculate the number of iterations $T$ required to reach $\epsilon$-accuracy, we note that the contribution of the gradient noise terms in (\ref{finalresult}) is given by
 \begin{equation}\label{var_red}
     \mathcal{O}\big(( \alpha +n\alpha^2  )\sigma_g^2/n\big) = \mathcal{O}\big( \alpha\sigma_g^2/n\big)\qquad \text{if }\alpha \leq 1/n
      \end{equation}
while the contribution  of the communication noise terms in (\ref{finalresult}) is given by:
    \begin{align}\label{arrstep}
     \mathcal{O}\Bigg( \left(\frac{(1+T\alpha)\gamma^2}{\alpha}+(\alpha^2 \ka^2 +\alpha^2T)n\ka \gamma^2 \right)\frac{\sigma_c^2}{n}\Bigg)
     &= \mathcal{\tilde{O}}\Bigg(\left( \frac{(1+T\alpha)\gamma^2}{\alpha}+ \frac{n \alpha^2}{\gamma} + n\alpha^2\gamma T\right)\frac{\sigma_c^2}{n}  \Bigg), 
    \end{align}
     where we used $\ka= \mathcal{\tilde{O}}(1/\gamma)$ and ignored the dependency on other problem parameters.  If we set $\gamma= \tilde{\mathcal{O}}(\alpha) $ such that (\ref{50}) is satisfied, the above bound further simplifies to 
    $$
   \mathcal{\tilde{O}}\Bigg(\left( \frac{(1+T\alpha)\gamma^2}{\alpha}+ \frac{n \alpha^2}{\gamma} + n\alpha^2\gamma T\right)\frac{\sigma_c^2}{n}  \Bigg)=  \mathcal{\tilde{O}}\left(\left((1+T\alpha)\alpha + n\alpha+n\alpha^3T\right)\frac{\sigma_c^2}{n}\right)
    $$
   For any given $\epsilon>0$, we can set $\alpha=\epsilon$ implying that $T = \mathcal{\tilde{O}}(\epsilon^{-1})$ iterations are required to achieve the specified $\epsilon$-accuracy.
\end{remark}

   \begin{remark} \textbf{($\mathbf{\sigma_c^2=0,\,\sigma_g^2=0}$ and $\mathbf{\sigma_c^2=0,\,\sigma_g^2>0}$):} In the absence of communication or gradient approximation errors ($\sigma_c^2=0,\,\sigma_g^2=0$), we can achieve the deterministic linear convergence rate of the gradient tracking algorithm \cite{DIG}. Referring to equation (\ref{finalresult}), we obtain the following inequality:
    $$
    \dd_{T} \leq(1-\alpha\mu/4 )^T\left(\dd_0 +  \frac{800(1+\ka^2) L}{n (1-\alpha\mu/4 )\mu}\| \Psi_{0}\|^2 \right)
    $$
The case $\sigma_c^2=0,\,\sigma_g^2>0$ considers stochastic decentralized optimization with no communication noise. For this scenario, with a constant $\alpha>0$, we have linear convergence to a neighbourhood of size $\mathcal{O}\big((\alpha^2 n + \alpha)\sigma_g^2/n\big)$ \cite{pudi}. A point to be remarked here is that \texttt{IC-GT} not only removes the data heterogeneity terms which arise in the convergence bound for \texttt{DGD} (cf. Table 1) but also makes sure that the variance scales linearly with the number of nodes provided $\alpha \leq 1/n$ (cf. \eqref{var_red}).

  \end{remark} 
\rb{
  \begin{remark} \textbf{(Consensus Error):} 
  We can establish convergence error bounds for the expected consensus error $\mathbb{E}[\|\Psi_k\|^2]$ by combining the results of Lemma \ref{r2} and Theorem \ref{thm1}. However, for brevity, we omit the explicit presentation of these results as they are of the same order as the results for $\Delta x^*_T$.
  \end{remark} 
}

\subsubsection*{Proof of Theorem \ref{thm1}}


Using (\ref{mat0}) and recalling that $\xx := \frac{(1_n\1_n^T) \otimes I_d}{n}\xb$,  the recursion for $\xx_k$ can be expressed as
\begin{align}\label{1-1}
\xx_{k+1} &= \bar{\textbf{v}}_k - \alpha \bar{\textbf{y}}_k \nonumber\\
&=\xx_{k}+ \gamma\bar{\bm{\epsilon}}_{k,c} -\alpha \bar{\textbf{y}}_k ,
\end{align}
where $\bar{\bm{\epsilon}}_{k,c}:=  \frac{1}{n} \big(1_n 1_n^T\otimes I_d \big)\hat{\textbf{Q}} \bm{\epsilon}_{k,c}$ and the last equality is due to $\bar{\textbf{v}}_k = \bar{\xb}_k+\gamma \bar{\bm{\epsilon}}_{k,c}$. Similarly, the recursion for $\bar{\mathbf{y}}_k:= \frac{1}{n} \big(1_n1_n^T\otimes I_d \big)y_k$ can be given as,
$$
\bar{\mathbf{y}}_{k} = \bar{\mathbf{y}}_{k-1} +   \frac{1}{n} \big(1_n1_n^T\otimes I_d \big) \big(\n \textbf{F} (\xb_{k},\bm{\xi}_{k}) - \n \textbf{F} (\xb_{k-1},\bm{\xi}_{k-1}) \big) + \gamma \bar{\bm{\epsilon}}_{k-1,c}.
$$
Taking telescopic sum from $0$ to $k$ leads to the following recursion: 
\begin{equation}\label{eq:bary}
\bar{\mathbf{y}}_{k} =  \frac{1}{n} \big(1_n1_n^T\otimes I_d \big)  \n \textbf{F} (\xb_{k},\bm{\xi}_{k}) +\gamma  \sum_{j=1}^{k}\bar{\bm{\epsilon}}_{j-1,c}
\end{equation}
since $\bar{\mathbf{y}}_0 =\frac{1}{n} \big(1_n1_n^T\otimes I_d \big)  \n \textbf{F}(\xb_{0},\bm{\xi}_{0}) $. Plugging \eqref{eq:bary} 
in \eqref{1-1}, we get,
\begin{align}\label{ll}
\xx_{k+1} &= \xx_{k}   + \gamma\bar{\bm{\epsilon}}_{k,c} -   \frac{\alpha}{n} \big(1_n1_n^T\otimes I_d \big) \n \textbf{F} (\xb_{k},\bm{\xi}_{k})   - \gamma \alpha \sum_{j=0}^{k-1}\bar{\bm{\epsilon}}_{j,c}\nonumber\\
&= \bar{\textbf{x}}_{k}-   \frac{\alpha}{n} \big(1_n1_n^T\otimes I_d \big) \n \textbf{f} (\xb_{k})   + \gamma\left( \bar{\bm{\epsilon}}_{k,c} - \alpha \sum_{j=0}^{k-1}\bar{\bm{\epsilon}}_{j,c}\right)+ \frac{\alpha}{n} \big(1_n1_n^T\otimes I_d \big)\Big( \n \textbf{f}(\xb_{k}) -  \n\textbf{ F} (\xb_{k},\bm{\xi}_k)\Big)\nonumber\\
&= \bar{\textbf{x}}_{k}-   \frac{\alpha}{n} \big(1_n1_n^T\otimes I_d \big) \n \textbf{f} (\xb_{k})  + \underbrace{\alpha \bm{\epsilon}_{k,g} + \gamma \bar{\bm{\epsilon}}_{k,c}}_{\delta_k}  - \alpha\gamma \sum_{j=0}^{k-1}\bar{\bm{\epsilon}}_{j,c}
\end{align}
where $\bm{\epsilon}_{k,g} $ is defined to be $\bm{\epsilon}_{k,g} := \frac{1}{n} \big(1_n1_n^T\otimes I_d \big)\Big( \n \textbf{f}( \xb_{k}) -  \n \textbf{F} (\xb_{k},\bm{\xi}_k)\Big)$ with $\EX\, [\bm{\epsilon}_{k,g}|\xb_k]  = 0 $ and $\EX [\|\bm{\epsilon}_{k,g} \|^2]\leq \sigma_g^2$ from Assumption \ref{asmp4}. Now, let $\mathcal{F}_k\defeq 
\sigma(\xb_{0}, \bm{\xi}_{0}, \bm{\epsilon}_{0,c}, \cdots, \bm{\xi}_{k-1}, \bm{\epsilon}_{k-1,c})$
be the sigma algebra generated by the random variables up to iteration $k$. Then, for any constant $\beta>0$, we have,
\begin{align}\label{altnew}
&\EX [ \|\bar{\mathbf{x}}_{k+1}-\mathbf{x}^*\|^2|\mathcal{F}_k] \nonumber \\
&~~~~~\leq (1+\beta)\EX [\|\bar{\mathbf{x}}_{k}-   \frac{\alpha}{n} \big(1_n1_n^T\otimes I_d \big) \n \textbf{f} (\xb_{k})   -\mathbf{x}^* + \delta_k\|^2 |\mathcal{F}_k]  + (1+\beta^{-1})\alpha^2 \gamma^2\EX \left[ \left\|\sum_{j=0}^{k-1}\bar{\bm{\epsilon}}_{j,c} \right\|^2\Big|\mathcal{F}_k\right]\nonumber\\
 &~~~~~=  (1+\beta) \|\bar{\mathbf{x}}_{k}-   \frac{\alpha}{n} \big(1_n1_n^T\otimes I_d \big)\n \textbf{f} (\xb_{k})   -\mathbf{x}^*\|^2  \nonumber\\
 &~~~~~\qquad + (1+\beta)\EX[ \|\delta_k\|^2 |\mathcal{F}_k]+ (1+\beta^{-1})\alpha^2 \gamma^2\EX \left[ \left\|\sum_{j=0}^{k-1}\bar{\bm{\epsilon}}_{j,c} \right\|^2\Big|\mathcal{F}_k\right],
\end{align}
where the equality is due to $\EX [\delta_k| \mathcal{F}_k] = 0$ from Assumption \ref{asmp3}. 
From Assumptions \ref{asmp2} and \ref{asmp4}, we have, 
\begin{equation}\label{del1}
      \EX [\|\delta_k\|^2] = \EX[\left\|\alpha \bm{\epsilon}_{k,g} + \gamma \bar{\bm{\epsilon}}_{k,c}\right\|^2]\leq\alpha^2 \sigma_g^2 + \gamma^2\sigma^2_{c}, 
\end{equation}
where we have used $\EX [\langle \bm{\epsilon}_{k,g}, \bar{\bm{\epsilon}}_{k,c}\rangle] = 0$. 
Furthermore, we have,
\begin{align}\label{epsbound}
  \EX \left[\left\| \sum_{j=0}^{k-1}\bar{\bm{\epsilon}}_{j,c} \right\|^2\right] &=  \EX \left[ \sum_{j=0}^{k-1}\left\|\bar{\bm{\epsilon}}_{j,c} \right\|^2\right]  + \sum_{1\leq p,p' \leq k-1}\EX \left[\langle \bar{\bm{\epsilon}}_{p,c} ,\bar{\bm{\epsilon}}_{p'
  ,c} \rangle \right]  
  \leq \sum_{j=0}^{k-1} \sigma_c^2 =  k\sigma_c^2,
\end{align}
where we use
$ \EX [\langle \bar{\bm{\epsilon}}_{p,c}, \bar{\bm{\epsilon}}_{p',c}\rangle] =  \EX [ \EX  \left[\langle \bar{\bm{\epsilon}}_{p,c}, \bar{\bm{\epsilon}}_{p',c}\rangle|\mathcal{F}_{p'}\right]]=0$ for $p<p'$.
Taking full expectations in (\ref{altnew}), it then follows that,
\begin{multline*}
 \EX  [\|\bar{\mathbf{x}}_{k+1}-\mathbf{x}^*\|^2 ]\leq (1+\beta)\EX \left[\left\|\bar{\mathbf{x}}_{k}-   \frac{\alpha}{n} \big(1_n1_n^T\otimes I_d \big)\n \textbf{f} (\xb_{k})  -\mathbf{x}^* \right\|^2\right]\\ + \big( (1+\beta)(\gamma^2\sigma_c^2 + \alpha^2\sigma_g^2) +k(1+\beta^{-1})\alpha^2\gamma^2\sigma^2_{c} \big).
\end{multline*}
where we used (\ref{del1}) to get the inequality.
We note that since $ \|\bar{\mathbf{x}}_{k+1}-\mathbf{x}^*\|^2 = \big\|\frac{(1\1^T) \otimes I_d}{n} (\textbf{x}_{k+1} - \textbf{x}^*)\big\|^2=n\|\bar{x}_{k+1}-x^*\|^2$, the above inequality leads to,
 \begin{multline}\label{ll0}
 \EX \left[\|\bar{x}_{k+1}-x^*\|^2\right] \leq (1+\beta) \EX  \left[ \left\|\bar{x}_{k}- \frac{\alpha }{n} \sum_{i=1}^n \n f_i (x_{i,k}) -x^* \right\|^2 \right]\\ + n^{-1}\big( (1+\beta)(\gamma^2\sigma_c^2 + \alpha^2\sigma_g^2) +k(1+\beta^{-1})\alpha^2\gamma^2\sigma^2_{c} \big).
\end{multline}
Considering the first term on the right hand side of \eqref{ll0}, we have,
\begin{align}\label{ll2}
    \Big\| \x_{k} - \frac{\alpha}{n} \sum_{i=1}^n\n f_i (x_{i,k}) -x^* \Big\|^2 &=  \| \x_k -x^* \|^2 -\frac{2\alpha}{n} \left\langle   \sum_{i=1}^n\n f_i (x_{i,k}), \x_{k}-x^* \right\rangle +\alpha^2 \Big\|\frac{1}{n}  \sum_{i=1}^n\n f_i (x_{i,k}) \Big\|^2.
\end{align}
The second term on the right hand side of \eqref{ll2} can be bounded as 
\begin{align}\label{bd2-}
     \langle   \sum_{i=1}^n\n f_i (x_{i,k}), \bar{x}_{k}-x^* \rangle &=  \langle   \sum_{i=1}^n\n f_i (x_{i,k}), \bar{x}_{k}-x_{i,k} \rangle +  \langle   \sum_{i=1}^n\n f_i (x_{i,k}), x_{i,k}-x^* \rangle \nonumber \\
     &\geq\sum_{i=1}^n \Big[ f_i(\bar{x}_k) -f_i(x_{i,k})- \frac{L}{2}\|\bar{x}_k -x_{i,k}\|^2 + f_i(x_{i,k}) - f_i(x^*)+ \frac{\mu}{2}\|x_{i,k}-x^*\|^2\Big] \nonumber\\
     &\geq\sum_{i=1}^n \Big[ f_i(\bar{x}_k) -f_i(x^*) -  \frac{L+\mu}{2} \|\bar{x}_k -x_{i,k}\|^2 +\frac{\mu }{4}\|\bar{x}_k -x^*\|^2\Big],
\end{align}
where the second inequality is due to Assumption~\ref{asmp3} and the last inequality is due to the inequality $ \|\bar{x}_k -x^*\|^2 \leq 2 \|\bar{x}_k-x_{i,k}\|^2 + 2\|x_{i,k} -x^*\|^2$. The last term on the right hand side of \eqref{ll2} can be bounded as
\begin{align}\label{bd1-}
     \Big\|\frac{1}{n}  \sum_{i=1}^n\n f_i (x_{i,k}) \Big\|^2&=  \Big\|\frac{1}{n}  \sum_{i=1}^n\n f_i (x_{i,k})-\frac{1}{n}  \sum_{i=1}^n\n f_i (\bar{x}_{k})+\frac{1}{n}  \sum_{i=1}^n\n f_i (\bar{x}_{k})-\frac{1}{n}  \sum_{i=1}^n\n f_i (x^*) \Big\|^2\nonumber\\
     &\leq  \frac{2L^2}{n} \sum_{i=1}^n \| \bar{x}_k-x_{i,k}\|^2+ \frac{4L}{n} \sum_{i=1}^n (f_i(\bar{x}_k) -f_i(x^*)),
\end{align}
where in the second summation, we used the fact that $\|\n f_i (\bar{x}_{k})  - \n f_i (x^*) \|^2 \leq 2L (f_i(\bar{x}_k) -f_i(x^*))$ by Assumption \ref{asmp3} \cite[Theorem 2.1.5]{yuri}. Using (\ref{bd2-}) and (\ref{bd1-}) in (\ref{ll2}) along with $\alpha<1/4L$, we have, 
\begin{align}\label{ll01}
    \left\| \bar{x}_{k+1} -  x^*\right\|^2 &\leq (1-\alpha\mu/2) \| \bar{x}_k -x^* \|^2  - \frac{\alpha}{n}\big(\sum_{i=1}^n f_i(\bar{x}_k) -f_i(x^*)\big) +\frac{(3L/2+\mu)\alpha}{n} \sum_{i=1}^n \|\bar{x}_k-x_{i,k}\|^2\nonumber\\
    &\leq (1-\alpha\mu/2) \| \bar{x}_k -x^* \|^2  - \alpha \left( \textbf{f}(\bar{\textbf{x}}_k) -\textbf{f}(\textbf{x}^*)\right) +\rb{\frac{5\alpha L}{2n}} \|\Psi_k\|^2,
\end{align}
where the last inequality is due to $\|\Delta \textbf{x}_k\|^2 \leq \|\Psi_k \|^2$.
Using (\ref{ll01}) in (\ref{ll0}), we get,
\begin{align*}
\EX[\|\bar{x}_{k+1}  -x^*\|^2]   &\leq (1+\beta)\Bigg\{(1-\alpha\mu/2) \EX [\| \x_k -x^* \|^2]  -  \alpha \left( \EX[\textbf{f}(\bar{\textbf{x}}_k)] -\textbf{f}(\textbf{x}^*)\right) \\ &\qquad +\rb{\frac{5\alpha L}{2n}} \EX[\|\Psi_k\|^2]\Bigg\} +\frac{1}{n}\big( (1+\beta)(\gamma^2\sigma_c^2 + \alpha^2\sigma_g^2) +k(1+\beta^{-1})\alpha^2\gamma^2\sigma^2_{c} \big).
\end{align*}
Set $\beta\defeq \frac{\alpha \mu}{4(1-\frac{\alpha \mu}{2})}$. We note that $(1+\beta^{-1}) \leq 4/\alpha \mu$ and $(1+\beta)= \frac{(1-\alpha\mu/4)}{(1-\alpha\mu/2)}$ with $1\leq (1+\beta)\leq 2$. Then, we have,
\begin{align}\label{r1}
\EX\left[\|\bar{x}_{k+1}  -x^*\|^2\right]   &\leq \Big\{(1-\alpha\mu/4) \EX \left[\| \x_k -x^* \|^2\right]  - \alpha \left( \EX[ \textbf{f}(\bar{\textbf{x}}_k)] -\textbf{f}(\textbf{x}^*)\right) \nonumber\\ &\qquad +\frac{5\alpha L}{n} \EX [\|\Psi_k\|^2]\Big\} +\frac{1}{n}\big(  (2+4\mu^{-1}k\alpha)\gamma^2 \sigma^2_{c} +2\alpha^2 \sigma_g^2\big). 
\end{align}
Multiplying both sides of \ref{r1} by $w_{k+1}\defeq (1 -\alpha \mu/4)^{-(k+1)}$, we have,
\begin{align}
w_{k+1}\dd_{k+1}&\leq 
(1-\alpha\mu/4) w_{k+1}\dd_{k} +\frac{5\alpha L}{n} w_{k+1}\EX [\|\Psi_k\|^2]\nonumber\\
&\qquad- \alpha w_{k+1}(\EX[ \textbf{f}(\mathbf{\x}_k)]-\textbf{f}(\mathbf{x}^*)) +\frac{w_{k+1}}{n}\big( \gamma^2 (2+4k\mu^{-1}\alpha)\sigma^2_{c} +2\alpha^2 \sigma_g^2\big) \nonumber.
\end{align}
Rearranging the terms, we get,
\begin{align}
 0  &\leq w_{k} \dd_k-w_{k+1}\dd_{k+1}  +\frac{5\alpha L}{n} w_{k+1} \EX [\|\Psi_k\|^2 ]\nonumber\\
&\qquad - \alpha w_{k+1}(\EX\, [\textbf{f}(\mathbf{\x}_k) ]-\textbf{f}(\mathbf{x}^*)) +\frac{w_{k+1}}{n}\big( \gamma^2 (2+4k\mu^{-1}\alpha)\sigma^2_{c} +2\alpha^2 \sigma_g^2\big) \nonumber \nonumber.
\end{align}
Summing the above inequality from $k=0$ to $T-1$, we get,
\begin{align}\label{intm22}
w_T\dd_{T} &\leq w_{0} \dd_0 +\frac{1}{n}\big( \gamma^2 (2+4T\mu^{-1}\alpha)\sigma^2_{c} +2\alpha^2 \sigma_g^2\big) \sum_{k=0}^{T-1}w_{k+1} + \frac{5\alpha L}{n}\sum_{k=0}^{T-1} w_{k+1}\EX[\|\Psi_k\|^2] \nonumber \\
&\qquad - \alpha \sum_{k=0}^{T-1} w_{k+1}(\EX [f(\mathbf{\x}_k)] -f(\mathbf{x}^*)).
\end{align}
We note that we can write the relations (\ref{main_rec})-(\ref{main_rec_2}) in Lemma \ref{r2} in the form of (\ref{rel}) 
with 
\begin{align}
   & b:=576 \alpha^2 L^2  \ka^2\qquad  c := 1344 \alpha^2\ka \nonumber \\
    &r:=   64\gamma^2 \big(2 + \alpha^2 (\ka^2 +1/2) + \alpha^2T\big)n\sigma_{c}^2 \ka+196n(\ka+1)\alpha^2\sigma^2_{g}\label{dd}\\
   & e_k:= \EX \left[\|\n \textbf{f}(\bar{\xb}_{k})-\n \textbf{f}(\xb^*) \|^2 \right]\nonumber
\end{align}
and we have taken $\rho'=1/4$ which fixes $\ka$ in (\ref{kappabd}) according to the bound (\ref{kappaprop}) (cf. (\ref{rho'eq})). Note that since $\alpha <  \frac{\sqrt{\rho'}}{2\sqrt{576} L \ka} $, $b\leq \frac{\rho'}{4}=\frac{1}{16}$. Then, with $a_t \defeq\EX[ \| \Psi_{t}\|^2] $ in Lemma \ref{lem4}, we get,
\begin{align}\label{bss1}
\EX [ \| \Psi_{t}\|^2] &\leq 40(1+\tau^2) \left(1- \frac{3\rho}{4\ka}\right)^t \|  \Psi_{0}\|^2 + 60 c  \sum_{j=0}^{t-1}   \left(1- \frac{3\rho}{4\ka}\right)^{t-j}  \EX \left[ \|\n \textbf{f(}\bar{\xb}_{j})-\n \textbf{f}(\textbf{x}^*) \|^2\right] + 52 r
\end{align}
with $\rho'' = 2(1+\ta^2)\|\Psi_0\|^2$ and $\rho=1-2\rho'=1/2$. We next bound the summation $\sum_{k=0}^{T-1} w_{k+1} \EX \|\Psi_k\|^2$ in (\ref{intm22}). To do this, we multiply both sides  of (\ref{bss1}) by $w_{k+1}$ and sum from $t=0$ to $T-1$: 
\begin{align}\label{intm101}
\sum_{k=0}^{T-1} (1-\alpha\mu/4)^{-(k+1)}  \EX \|\Psi_k\|^2 &\leq   40(1+\ka^2)\|\Psi_0\|^2\sum_{k=0}^{T-1} (1-\alpha\mu/4)^{-(k+1)} \left(1- \frac{3\rho}{4\ka}\right)^k\nonumber \\
&\qquad +   60c\sum_{k=0}^{T-1}(1-\alpha\mu/4)^{-(k+1)} \sum_{j=0}^{k-1} \left(1- \frac{3\rho}{4\ka}\right)^{t-j} e_j + 52 r  W_{T-1},
\end{align}
where $W_{T-1} =  \sum_{k=0}^{T-1} w_{k+1}$.  From (\ref{alphacond}), we have,
 \begin{equation}\label{a-r}
     \alpha \mu/2 \leq 3\rho/4\ka \implies  \alpha \mu/2(1-\alpha\mu/8) \leq 3\rho/4\ka \implies 1-3\rho/4\ka \leq (1-\alpha \mu/4)^2.
 \end{equation}
 We use (\ref{a-r}) to bound the two summations on the right hand side of (\ref{intm101}) as follows:
 \begin{equation}\label{fbd1}
\sum_{k=0}^{T-1} (1-\alpha\mu/4)^{-(k+1)} \Big(1- \frac{3\rho}{4\ka}\Big)^k \leq \sum_{k=0}^{T-1} (1-\alpha\mu/4)^{k-1} \leq \frac{4w_1}{\alpha \mu}, 
\end{equation}
and 
\begin{align}
     \sum_{k=0}^{T-1}(1-\alpha\mu/4)^{-(k+1)} &\sum_{j=0}^{k-1} \Big(1- \frac{3\rho}{4\ka}\Big)^{k-j} e_j= \sum_{k=0}^{T-1}\sum_{j=0}^{k-1} (1-\alpha\mu/4)^{-(k+1)+j+1}  \Big(1- \frac{3\rho}{4\ka}\Big)^{k-j} w_{j+1} e_j \nonumber\\
     &~~~~~= \sum_{k=0}^{T-1}\sum_{j=0}^{k-1}  \Bigg(\frac{1- 3\rho/4\ka}{1-\alpha\mu/4}\Bigg)^{k-j} w_{j+1} e_j 
     \leq
     \sum_{k=0}^{T-1}\sum_{j=0}^{k-1}  (1-\alpha\mu/4)^{k-j} w_{j+1} e_j \nonumber \\
     &~~~~~\leq \sum_{k=0}^{T-1}  (1-\alpha\mu/4)^{k} \sum_{k=0}^{T-1} w_{k+1} e_k \leq \frac{4}{\mu \alpha}\sum_{k=0}^{T-1} w_{k+1} e_k \label{fbd2},
\end{align}
where the first equality is due to \eqref{a-r} and the second inequality is obtained using the relation $\sum_{k=0}^{T-1} \sum_{j=0}^{k-1}a_{k-j}b_j \leq\sum_{k=0}^{T-1} a_{k} \sum_{k=0}^{T-1} b_{k}$ for any two non-negative scalar sequences $a_k,\,b_k,\,k\geq 0$. Plugging the previous two bounds in (\ref{intm101}), we get,
\begin{align*}
\sum_{k=0}^{T-1} w_{k+1}  \EX[ \|\Psi_k\|^2] &\leq  \frac{ 160w_1(1+\ka^2) \|\Psi_0\|^2 }{\alpha\mu }+\frac{240 n cL}{\mu \alpha}\sum_{k=0}^{T-1}w_{k+1} \left(\EX[ \textbf{f}(\bar{\xb}_k)] -\textbf{f}(\textbf{x}^*)\right) +  52r W_{T-1},
\end{align*}
where we have additionally used the fact that $\|e_k\|^2 = \EX \left[ \|\n \textbf{f}(\mathbf{\x}_k))-\n \textbf{f}(\textbf{x}^*) \|^2 \right]\leq 2nL( \EX [\textbf{f}(\mathbf{\x}_k))] -\mathbf{f}(\textbf{x}^*)) $ from Assumption \ref{asmp3} \cite[Theorem 2.1.5]{yuri}. Finally, using the above bound in (\ref{intm22}), we get,
\begin{align*}
  w_T\dd_{T} &\leq w_{0} \dd_0 +\frac{1}{n}\big( \gamma^2 (2+4T\mu^{-1}\alpha)\sigma^2_{c} +2\alpha^2 \sigma_g^2\big) W_{T-1} \\
&\quad +\frac{5\alpha L}{n}\Bigg( \frac{160w_1(1+\ka^2)\|\Psi_0\|^2 }{\mu \alpha}+\frac{240ncL}{\mu \alpha}\sum_{k=0}^{T-1}w_{k+1} \left(\EX [\textbf{f}(\mathbf{\x}_k)) ]-\mathbf{f}(\textbf{x}^*)\right) +  52r W_{T-1}\Bigg) \\
 & \quad - \alpha \sum_{k=0}^{T-1}w_{k+1}(\EX [\,\textbf{f}(\mathbf{\x}_k)] -\mathbf{f}(\textbf{x}^*)).
\end{align*}
 Rearranging the terms in the above inequality and recalling that $c=1344 \alpha^2\ka$, we get, 
\begin{align*}
\dd_{T} &\leq \frac{1}{w_T} \Bigg\{\dd_0 + \frac{800w_1(1+\ka^2) L}{n \mu } \| \Psi_{0}\|^2 \Bigg\}+\frac{2}{\mu\alpha}\Bigg\{  \frac{\gamma^2 (2+4T\mu^{-1}\alpha)\sigma^2_{c}}{n} +\frac{2\alpha^2\sigma_g^2}{n}  \Bigg\}  \\
&\quad+\frac{ 520  L r}{ n\mu }  +\alpha\underbrace{\Bigg\{\frac{1612800\ka L^2}{\mu }\alpha-1\Bigg\}}_{\leq 0} \sum_{k=0}^{T-1}\frac{w_{k+1}}{w_T} \underbrace{(\EX [\mathbf{f}(\mathbf{\x}_k)] -f(\mathbf{x}^*))}_{\geq 0},
\end{align*}
where we used $ W_{T-1}/w_T \leq2/  \mu \alpha$. The last term on the right had side is less than zero due to the condition on $\alpha$ (see \eqref{alphacond}). Plugging the value of $r$ from (\ref{dd}) in the above inequality completes the proof.\qed 


\section{Numerical Experiments}\label{sec:num}
In this section, we present an empirical evaluation of the performance of \texttt{IC-GT} through two sets of numerical experiments. 
The first set focuses on logistic regression on the MNIST dataset, while the second set 
explores the effect of different 
noise variances in a deep learning setting. All experiments were implemented using PyTorch, with a dedicated CPU core functioning as a node.

\subsection*{Logistic regression}
We first consider $\ell_2$ regularized logistic regression problems of the form,
  \begin{equation}
       \min_{x\in \mathbb{R}^d} \, \left\{L(x;y,z) :=-\frac{1}{m} \sum_{i=1}^m\Big\{ z^i \log \varphi(x^Ty^i) + (1-z^i)\log (1-\varphi(x^Ty^i))\Big\} + \frac{\lambda}{2}\|x\|^2\right\}, 
    \end{equation}
where  $x \in \mathbb{R}^d$ denote the learnable model parameters, 
$\{y^i,z^i\}_{i=1}^m$ denote the set of $m$ data points, $\varphi(\cdot)$ denotes the sigmoid function, and $\lambda>0$  is the regularization parameter. We use the \texttt{MNIST} dataset which consists of 60,000,  28$\times$28 pixel grayscale images of handwritten single digits between 0 and 9. The data is partitioned in a disjoint manner amongst the nodes by assigning each node $10^3$ data samples independently. 
\begin{figure}[H]
    \centering
    \begin{minipage}[b]{0.4\textwidth}
        \centering
         \includegraphics[width=\textwidth]{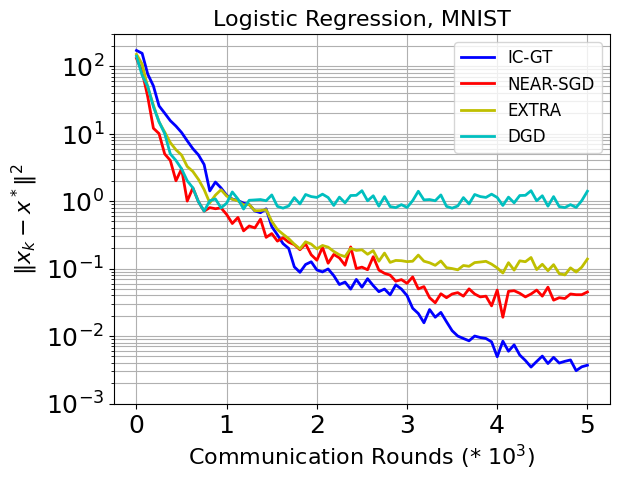}
         \subcaption{}
        \label{fig:plot1}
    \end{minipage}
    \begin{minipage}[b]{0.4\textwidth}
        \centering
        \includegraphics[width=\textwidth]{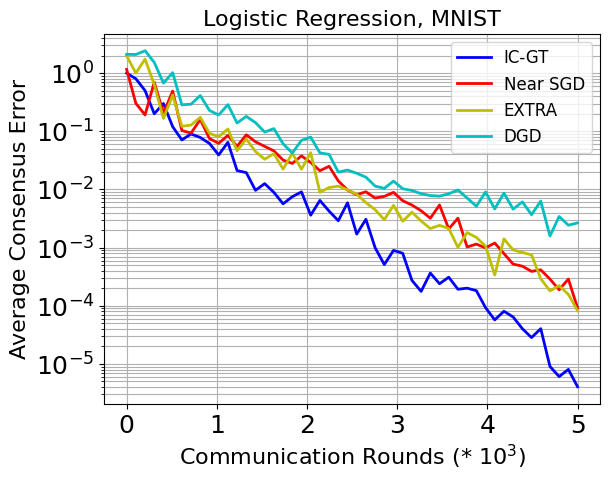}
         \subcaption{}
        \label{fig:plot2}
    \end{minipage}
    \hfill
    \begin{minipage}[b]{0.4\textwidth}
        \centering
        \includegraphics[width=\textwidth]{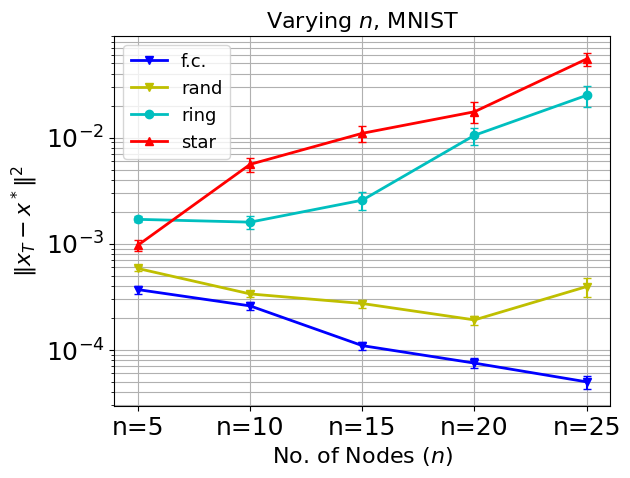}
         \subcaption{}
        \label{fig:plot3}
    \end{minipage}
    \label{fig:four-side-by-side}
      \centering
    \begin{minipage}[b]{0.41\textwidth}
        \centering
        \includegraphics[width=\textwidth]{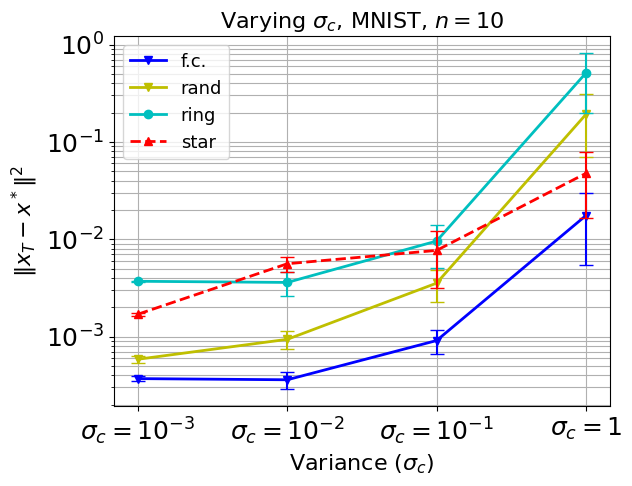}
        \subcaption{}
        \label{fig:plot4}
    \end{minipage}
     \caption{\footnotesize{\textbf{(a)-(b)}  Optimality Error, Average consensus error vs. communication rounds for \ty{MNIST} dataset with star topology $(n=10)$. \textbf{(c)} Final Optimal error $\|x_T-x^*\|^2,\,T=5\times 10^3$ for $n\in \{5,10,15,20,25\}$ for different topologies. \textbf{(d)} Final Optimal error $\|x_T-x^*\|^2$ for $\sigma_c \in \{10^{-3},10^{-2},10^{-1},1\}$ }}
    \label{fg1}    
\end{figure}
To simulate the inexact communication setting, we incorporate zero-mean Gaussian noise with a variance of $\sigma_c^2$ into the transmitted model estimates independently. We adopt a star topology with $n=10$ for the communication structure. In evaluating the performance, we employ the $\ell_2$ distance between the averaged variable $\x_k$ and the optimal point $x^*$. The optimal point $x^*$ is computed using the L-BFGS algorithm from the SciPy library in Python.
We also include the average consensus error as a performance metric, which is computed as $\frac{1}{|\mathcal{E}|}\sum_{(i,j)\in \mathcal{E}} \|x_i-x_j\|^2$, where $\mathcal{E}$ represents the edge set. We compare our proposed algorithm (\ty{IC-GT}) with several baselines, including the \ty{NEAR-SGD} algorithm from \cite{do4}, the \ty{EXTRA} algorithm proposed in \cite{EXTRA}, and the gradient tracking method \cite{sgt1}. Additionally, we include the performance of the \ty{DGD} algorithm for comparison purposes.

In our experiments, we set the batch size to $32$ and tune the step size $\alpha$ using a grid-search over the range $\alpha \in [10^{-4}, 1]$ to obtain the best performance for all the algorithms. The total number of communication rounds is set to $T=5\times 10^3$. For \ty{IC-GT}, we set the attenuation noise parameter $\gamma$ to $\gamma = \alpha \times \log T$.
The performance results are reported in Figure \ref{fg1}(a)-(b). From the plots, it is evident that \ty{IC-GT} outperforms all the other algorithms in terms of both the optimality error and the consensus error.

To assess the scalability of \ty{IC-GT} and examine the impact of graph connectivity on its convergence accuracy, we conducted experiments with varying network sizes, specifically $n \in \{5, 10, 15, 20, 25\}$. We kept the noise variance fixed at $\sigma_c^2=0.01$ for the following graph topologies: (i) Fully connected (f.c.), (ii) Erdős-Rényi graph with an edge probability of 0.5 (rand), (iii) Ring topology, and (iv) Star topology. From Figure \ref{fg1}(c), we observe that as the graph connectivity deteriorates, the final performance of \ty{IC-GT} also deteriorates. In the case of a fully connected graph, there is an improvement in performance with an increasing number of nodes due to a decrease in gradient variance resulting from an increased effective mini-batch size. Finally, we also investigate the effect of varying $\sigma_c^2$ on the performance of \ty{IC-GT}, as depicted in Figure \ref{fg1}(d).

\subsection*{Neural network based experiments}
In this subsection, we investigate a deep learning scenario that involves random compressed communication using probabilistic quantization (see \eqref{randquant}). We assume a star-based topology with $n=10$ for both the \ty{MNIST} and \ty{CIFAR} datasets. For the \ty{MNIST} dataset, we utilize a learning model with a total of $8.4$K parameters. This model comprises two convolution layers, the first with $250$ parameters and the second with $5$K parameters, followed by a fully connected layer with $3.2$K parameters. For the \ty{CIFAR-10} dataset, we adopt the standard \ty{LENET} architecture, which consists of three convolution layers and two fully connected layers. This architecture has a total of $0.54$M parameters. The configuration of the max-pooling and batch normalization layers follows the standard settings used in \ty{LENET} models.


We compare \ty{IC-GT} with two other strategies commonly employed to address noise in an inexact communication setting. The first strategy involves utilizing a decreasing noise variance policy, where the variance decreases as the number of communication rounds progresses. In this approach, we employ \ty{GT} with quantization and adjust the quantization levels to become finer as the rounds increase. Specifically, in the case of \eqref{randquant}, we uniformly increase the parameter $\Delta_p$ from $\Delta_p=1$ to $\Delta_p=5 \times 10^3$ as the rounds progress. This results in higher levels of noise variance in the initial rounds and lower levels in the final rounds. The second strategy maintains a uniform quantization level of $\Delta_p=10^2$ throughout all communication rounds, leading to a fixed noise variance. We employ the same quantization level of $\Delta_p=10^2$ for \ty{IC-GT}.

\begin{figure}[htp]
    \centering
    \begin{minipage}[b]{0.4\textwidth}
        \centering
        \includegraphics[width=\textwidth]{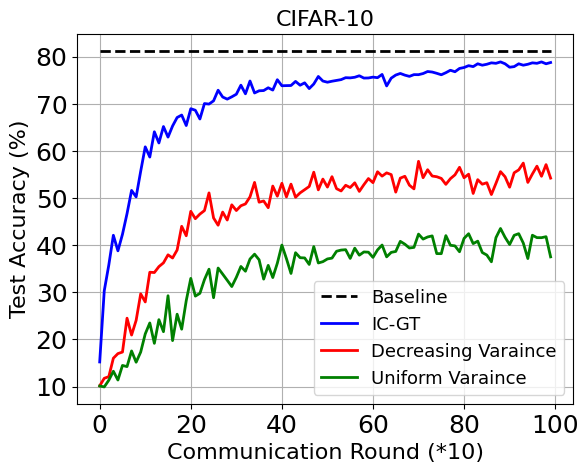}
          \subcaption{}
        \label{fig:plot5}
    \end{minipage}
    \begin{minipage}[b]{0.41\textwidth}
        \centering
        \includegraphics[width=\textwidth]{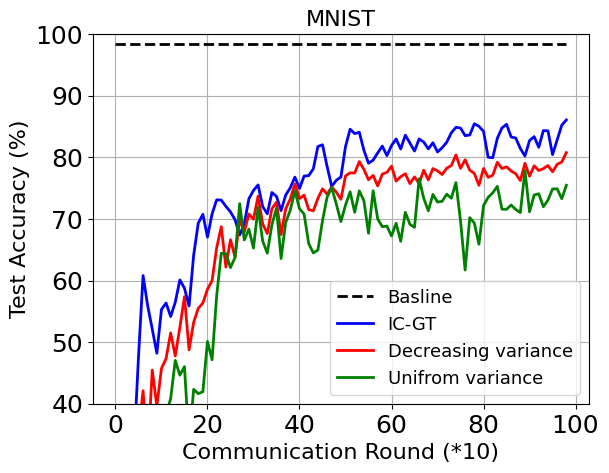}
        \subcaption{ }
        \label{fig:plot6}
    \end{minipage}
     \caption{\footnotesize{\textbf{(a)} Test Accuracy vs. communication rounds for CIFAR-10 dataset with star topology ($n=10$).  \textbf{(b)} Test Accuracy vs. communication rounds for MNIST dataset with star topology ($n=10$).}}
    \label{fg2}    
\end{figure}

The results of the comparison have been plotted in Figure~\ref{fg2}(a)-(b). In both plots, the baseline represents the highest achievable accuracy that can be obtained in a centralized setting using the models employed. From the plots, we observe that for both the \ty{CIFAR-10} and \ty{MNIST} datasets, the performance of \ty{IC-GT} is the closest to the baseline. The performance difference between \ty{IC-GT} and the baseline appears to be more pronounced in the case of \ty{CIFAR-10} compared to \ty{MNIST}.

\section{Final Remarks}
In this paper, we proposed a gradient tracking based algorithm for decentralized optimization in an inexact communication scenario. We established theoretical convergence guarantees and analyzed the impact of communication and gradient noise on performance. Our algorithm effectively mitigates the impact of communication noise and data heterogeneity,  and achieves optimal iteration complexity for strongly convex, stochastic smooth functions. Experimental results on logistic regression and neural networks demonstrated the superiority of the proposed algorithm over existing methods. As future work, the algorithm can be extended to other settings, such as directed graphs and asynchronous updates, and incorporate variance reduction techniques to enhance convergence rate.
\section*{Appendix I: Proof of Lemma \ref{lem0}}\label{sec:apndxI}
\begin{proof}
	From (\ref{mat0}), we have,
	\begin{align}\label{mat}
		\vb_{k}&= (\II-\gamma \Q)\xb_{k} + \gamma \hat{\textbf{Q}}\epsilon_{k,c}.
	\end{align}
	Multiplying both sides of (\ref{mat}) by $\frac{1}{n}1_n1_n^T\otimes I_d$, we get,
	\begin{align}\label{mat1}
		\Bar{\textbf{v}}_k &= \Bar{\textbf{x}}_{k} +  \frac{\gamma}{n} \big(1_n1_n^T\otimes I_d \big)\hat{\textbf{Q}} \epsilon_{k,c} \nonumber\\
		& = (\II-\gamma \Q) \bar{\textbf{x}}_k  + \frac{\gamma}{n} \big(1_n 1_n^T\otimes I_d \big)\hat{\textbf{Q}} \epsilon_{k,c}.
	\end{align}
	where we used $\frac{1}{n}1_n^T\Q =0$ to get the first inequality and $(\II-\gamma \Q)\bar{\textbf{x}}_k  = \bar{\textbf{x}}_k $ to get the last inequality. Subtracting (\ref{mat1}) from  (\ref{mat}) and adding $- \frac{1}{n}1_n1_n^T \Delta \xb_{k}$, we get,
	\begin{equation}\label{mat3}
		\Delta \vb_k = (\III-\gamma \Q ) \Delta \xb_{k} +  \gamma \tilde{\textbf{Q}}\epsilon_{k,c}. 
	\end{equation}
	From (\ref{matI}), the expression for $\Delta \xb_k$ can be written as,
	\begin{equation}\label{mat4}
		\Delta \xb_{k} = \Delta \vb_{k-1}  - \alpha \Delta \yb_{k-1}.
	\end{equation}
	Substituting for $\Delta \xb_k$ in (\ref{mat3}) using (\ref{mat4}) yields the following
	recursive relation for $\Delta \vb_k$ in terms of $\Delta \vb_{k-1}$ and $\Delta \yb_{k-1}$:
	\begin{align*}
		\Delta \vb_k =(\III-\gamma \Q ) \Delta \vb_{k-1}  - \alpha (\III-\gamma \Q ) \Delta \yb_{k-1} +  \gamma \tilde{\mathbf{Q}} \epsilon_{k,c}
	\end{align*}
	Next, the recursive relation for $\Delta \xb_k$ in terms of $\Delta \xb_{k-1}$ and $\Delta \yb_{k-1}$ is obtained by substituting  for $\Delta \vb_{k-1}$ in (\ref{mat4}) using \eqref{mat3}. That is,
	\begin{align*}
		\Delta \xb_k =(\III-\gamma \Q ) \Delta \xb_{k-1}  +  \gamma \tilde{\textbf{Q}} \epsilon_{k-1,c} - \alpha \Delta \yb_{k-1}. 
	\end{align*}
The recursive expression for $\Delta \yb_k$ can be obtained similarly using the expression for $\bar{\yb}_k$ and subtracting it from (\ref{matII}), concluding the proof. 
\end{proof}

\subsection*{Appendix II: Proof of Lemma \ref{lem3}}\label{sec:apndxII}
\begin{proof}
	Using mathematical induction, we can show that $\textbf{J}^\tau_\gamma$ for any $\ka\in \mathbb{N}$ is given as, 
	\begin{equation}\label{matproddef}
		\textbf{J}^\ka_\gamma =
		\begin{bmatrix}
			(\III-\gamma \Q )^\ka &0 &-\ka (\III-\gamma \Q )^\ka \\
			0 & (\III-\gamma \Q)^{\ka} &-\ka (\III-\gamma \Q )^{\left(\ka-1\right)} \\
			0 & 0&(\III-\gamma \Q)^{\ka} \
		\end{bmatrix}.
	\end{equation}
	Taking norms in \eqref{matproddef} and using triangular inequality, we get,
	\begin{equation}\label{basic}
		\|\textbf{J}_\gamma^\ka\| \leq \|(\III-\gamma \Q )^{\ta}\| +\ta\|(\III-\gamma \Q)^{(\ta-1)}\| +\ta\|(\III-\gamma \Q)^{\ta}\|.    
	\end{equation}
	We will next bound the terms on the right hand side of \eqref{basic}. 
Note that the smallest eigenvalue of the matrix $(\III-\gamma \Q)^\ka$ is zero and the remaining eigenvalues 
	are of the form $(1-\gamma (1-\lambda_i))^\ka$ for $i =2,\dots,n$, where $\lambda_i$ are the eigenvalues of $Q$ defined in Assumption~\ref{asmp1}. Therefore, 
	\begin{align}\label{2normbd}
		\|(\III-\gamma \Q )^\ka\| &= \max_{i=2,\dots,n} (1-\gamma(1-\lambda_i))^{\ka}\nonumber \\
		&= (1-\gamma(1-\lambda_2))^{\ka}.
	\end{align}
 From \eqref{kappaprop}, it follows that $	\ka \geq 2\left(1-\frac{\ln \sqrt{\delta}/4}{ \gamma (1-\lambda_2)}\right) > -\frac{\ln \sqrt{\delta}/2}{ \gamma (1-\lambda_2)}$. Substituting this inequality in \eqref{2normbd}, we get, 
	\begin{align}\label{5-2}
			\big(1-\gamma (1-\lambda_2)\big)^{\ka} &\leq \exp\big( -\ka\gamma  (1-\lambda_2) \big)\leq \frac{\sqrt{\delta}}{2}.
	\end{align}	
We next bound the second term in \eqref{basic}. For convenience, we define $\textbf{Q}_1\defeq\ka (\III-\gamma \Q )^{\ka-1} $. The smallest eigenvalue of $\textbf{Q}_1$ is zero and the remaining eigenvalues are of the form $\ka (1-\gamma(1-\lambda_i) )^{\ka-1},$ for $i=2,\dots,n$. Therefore,  
	\begin{align}\label{Amatrix_1}
		\|\textbf{Q}_1\| &\leq \ka (1-\gamma(1-\lambda_2) )^{\ka-1} 
		\leq \ka \exp (- (\ka-1)\gamma(1-\lambda_2) ). 
	\end{align}
	 Taking logarithm on both sides of \eqref{Amatrix_1} yields,
	 \begin{align}\label{Amatrix}
	 			\ln \|\textbf{Q}_1\| &\leq  \ln \ka - (\ka-1)\gamma(1-\lambda_2).
	 \end{align}	
	Now, consider $\frac{\ln \ka}{\ka}$  as a function of $\ka$ and observe that it is monotonically decreasing for any $\ka > \exp(1)$ since its first derivative $\frac{1 - \ln \ka}{\ka^2} < 0$. From \eqref{kappaprop}, we have $\tau\geq 16 \ln 4 > \exp(1)$ since $\gamma < 1/4$ and $\lambda_2 \in (-1,1)$.  For convenience, we define $\epsilon_{\gamma,\lambda_2} = \frac{\gamma (1 - \lambda_2)}{2} \in [0,1/4)$.  Therefore, from\eqref{kappaprop}, it follows that,  
	\begin{align}\label{tmp9}
		\frac{\ln \tau}{\tau} \leq \epsilon_{\gamma,\lambda_2} \frac{\ln \frac{4}{\epsilon_{\gamma,\lambda_2}}+ \ln \ln \frac{1}{\epsilon_{\gamma,\lambda_2}}}{4\ln 1/\epsilon_{\gamma,\lambda_2}}\leq \epsilon_{\gamma,\lambda_2} =  \frac{\gamma (1 - \lambda_2)}{2}.
	\end{align}	
Using (\ref{tmp9}) and $	\ka \geq 2\left(1-\frac{\ln \sqrt{\delta}/4}{ \gamma (1-\lambda_2)} \right)$ in  \eqref{Amatrix}, we get, 
\begin{align}
	\ln \|\textbf{Q}_1\| &\leq  \ln \ka - (\ka-1)\gamma(1-\lambda_2) \nonumber \\
	&\leq \frac{\ka\gamma(1 - \lambda_2)}{2} -  (\ka-1)\gamma(1-\lambda_2) \nonumber \\
	&=   \gamma(1 - \lambda_2)\left(1-\frac{\ka}{2}\right) \nonumber \\ &\leq \ln \sqrt{\delta}/4. \nonumber 
\end{align}
Therefore, 
\begin{align}\label{Bmatrix}
	 \|\textbf{Q}_1\| \leq \sqrt{\delta}/4.
\end{align}	
Finally, we bound the third term in \eqref{basic} as, 
\begin{align}\label{Cmatrix}
	\ta\|(\III-\gamma \Q)^{\ta}\| \leq \|\textbf{Q}_1\| \leq \sqrt{\delta}/4.
\end{align}	
Combining, \eqref{basic}, \eqref{5-2},\eqref{Bmatrix} and \eqref{Cmatrix}, we get, 
	$$\|\textbf{J}^\ka_\gamma\|^2 \leq  \left(\frac{\sqrt{\delta}}{2} + \frac{\sqrt{\delta}}{4} + \frac{\sqrt{\delta}}{4}\right)^2 = \delta.
$$.

\end{proof}

\section*{Appendix III: Proof of Lemma \ref{r2}}\label{sec:apndxIII}
\begin{proof}
We begin by iterating the relation \eqref{mm1} with $k=t+\tau$:
\begin{align}
\Psi_{t+\ka} &=  \textbf{J}_\gamma \Psi_{t+\ka-1} + \alpha \textbf{E}_{t+\ka-1} \nonumber\\
&= \textbf{J}^\ka_\gamma \Psi_t + \alpha \sum_{i=0}^{\ka-1} \textbf{J}^{\ka-i-1}_\gamma  \textbf{E}_{t+i}\label{e0}
\end{align}
We next consider $\textbf{E}_k$ whose definition is recalled here:
$$
\textbf{E}_{k-1} =  \underbrace{\frac{ \gamma}{\alpha}
\begin{bmatrix}
\tilde{\textbf{Q}}\bm{ \epsilon}_{k,c} \\
\tilde{\textbf{Q}} \bm{\epsilon_{k-1,c}} \\
\alpha \tilde{\textbf{Q}} \bm{\epsilon_{k-1,c} }
\end{bmatrix} }_{\textbf{E}^c_{k-1}}
 +
\underbrace{ 
\begin{bmatrix}
0\\
0\\
\III (\nabla \textbf{F}(\xb_{k},\bm{\xi}_{k} ) - \nabla \textbf{F}(\xb_{k-1},\bm{\xi}_{k-1} ))
\end{bmatrix}}_{\textbf{E}_{k-1}^g}, \quad \forall k \in \N.
$$
We note that 
\begin{equation}\label{decomp1}
\EX \left[ \left\| \sum_{i=0}^{\ka-1} \textbf{J}^{\ka-i-1}_\gamma \textbf{E}_{t+i} \right\|^2\right] \leq  2\EX \left[\left\|\sum_{i=0}^{\ka-1} \textbf{J}^{\ka-i-1}_\gamma  \textbf{E}^c_{t+i} \right\|^2\right] + 2\EX \left[\left\| \sum_{i=0}^{\ka-1} \textbf{J}^{\ka-i-1}_\gamma \textbf{E}^g_{t+i}\right\|^2\right]
\end{equation}
We first bound the first term on the right hand side of \eqref{decomp1}.
Using the expression for the matrix product $\textbf{J}_\gamma^{\tau-i-1}$ for any $0\leq i\leq\ka-1$ (cf. (\ref{matproddef})), we have,
\begin{equation}\label{normmat}
\textbf{J}^{\ka-i-1}_\gamma \textbf{ E}^c_{t+i} = \frac{ \gamma}{\alpha}
\begin{bmatrix}
(\III-\gamma \Q)^{\ka-i-1} \tilde{\textbf{Q}} \bm{\epsilon}_{t+i+1,c} -  \alpha (\ka-i-1)(\III-\gamma \Q)^{\ka-i-1}\tilde{\textbf{Q}} \bm{\epsilon}_{t+i,c} \\
(\III-\gamma \Q)^{\ka-i-1} \tilde{\textbf{Q}}- \alpha (\ka-i-1) (\III-\gamma \Q)^{\ka-i-2} \tilde{\textbf{Q}}\big)\bm{\epsilon}_{t+i,c}  \\
 \alpha (\III-\gamma \Q)^{\ka-i-1} \tilde{\textbf{Q}}\bm{\epsilon}_{t+i,c}
\end{bmatrix}. 
\end{equation}
Note that, using $\|Q\|,\,\|\text{diag}(Q)\|\leq 1$, we have,
\begin{equation}\label{bds1}
\|(\III-\gamma \Q)^{\ka-i-1}\|^2\leq 1, \qquad  \|\hat{\textbf{Q}}\| = \|(Q' - \text{diag}(Q'))\otimes I_d\| \leq  2, \quad \mbox{and}\quad \|\tilde{\textbf{Q}} \| \leq
\left\|\II - \frac{11^T}{n}\right\|\|\hat{\textbf{Q}}\|  \leq 2. 
\end{equation}
Taking norms in (\ref{normmat}) and using the bounds (\ref{bds1}), we get,
\begin{align}\label{tmp111}
\EX\left[ \|\textbf{J}^{\ka-i-1}_\gamma  \textbf{E}^c_{t+i}\|^2\right]  &\leq \frac{4\gamma^2}{\alpha^2}( 2(1+ \alpha^2(\ka-i-1)^2) +\alpha^2) \max\{  \EX\left[ \|\bm{\epsilon}_{t+i+1,c} \|^2\right], \EX \left[\|\bm{\epsilon}_{t+i,c} \|^2 \right]\}\nonumber\\
&\leq \frac{4\gamma^2}{\alpha^2}( 2+ 2\alpha^2\ka^2 +\alpha^2) \max\{  \EX \left[\|\bm{\epsilon}_{t+i+1,c} \|^2\right], \EX\left[\|\bm{\epsilon}_{t+i,c} \|^2 \right]\}\nonumber\\
&\leq \frac{8\gamma^2}{\alpha^2}( 1+ \alpha^2\left( \ka^2 +1/2\right)n\sigma_{c}^2,
\end{align}
where the last inequality is due to Assumption \ref{asmp2}. Hence,  
\begin{align}\label{bd1}
   \EX \left[\left\|\sum_{i=0}^{\ka-1} \textbf{J}^{\ka-i-1}_\gamma \textbf{ E}^c_{t+i}\right\|^2\right] &= 
   \sum_{i=0}^{\ka-1}  \EX \left[\left\|\textbf{J}^{\ka-i-1}_\gamma  \textbf{E}^c_{t+i}\right\|^2\right]
   \leq \frac{8\gamma^2}{\alpha^2} \big(1 + \alpha^2 (\ka^2 +1/2)\big)n\sigma_{c}^2 \ka,
\end{align}
where the equality is due to Assumption \ref{asmp2} and the fact that the cross terms of the form $\langle \bm{\epsilon}_{i,c},\bm{\epsilon}_{j,c}\rangle $ are all zero. That is, 
if we denote $\mathcal{F}_k\defeq 
\sigma(\xb_{0}, \bm{\xi}_{0}, \bm{\epsilon}_{0,c}, \cdots, \bm{\xi}_{k-1}, \bm{\epsilon}_{k-1,c})$ 
to be the sigma algebra generated by the random variables up to iteration $k$, we have for any $i,j$ with $i<j$, $ \EX [\langle\bm{\epsilon}_{i,c}, \bm{\epsilon}_{j,c}\rangle] =  \EX [ \EX  \left[\langle\bm{\epsilon}_{i,c}, \bm{\epsilon}_{j,c}\rangle|\mathcal{F}_{j}\right]]=0$.

Next we consider $\textbf{E}^g_k$ to bound the second summation in (\ref{decomp1}). Let $\g_{k} \defeq  \n \textbf{F}(\xb_{k},\bm{\xi}_{k})-\n \textbf{f}(\xb_k)$ and $\ddd_k = \n \textbf{f}(\textbf{x}_k) - \n \textbf{f}(\textbf{x}^*)$. We note from Assumption \ref{asmp4}, $\g_k$ is a zero mean vector given $\xb_k$ with variance $n \sigma_g^2$. Using $\bar{\mathbf{Q}} 
\defeq(\III-\gamma \Q)$ 
and the expression for the matrix product $\textbf{J}_\gamma^{\ka-i-1}$ (cf. \ref{matproddef}), we have,
\begin{align}\label{sixterms}
 \EX &\left[\left\|\sum_{i=0}^{\ka-1} \textbf{J}^{\ka-i-1}_\gamma  \textbf{E}^g_{t+i}\right\|^2\right] \nonumber \\&= \EX \left[\left\|\sum_{i=0}^{\ka-1} (\ka-i - 1) \bar{\mathbf{Q}} ^{\ka-i-1}\III( \n \textbf{F}(\xb_{t+i+1},\bm{\xi}_{t+i+1})-\n \textbf{F}(\xb_{t+i},\bm{\xi}_{t+i}) ) \right\|^2\right] \nonumber \\
 &~~~~~+\EX\left[ \left\|\sum_{i=0}^{\ka-1}  (\ka-i - 1) \bar{\mathbf{Q}}^{\ka-i-2}\III(\n \textbf{F}(\xb_{t+i+1},\bm{\xi}_{t+i+1})-\n \textbf{F}(\xb_{t+i},\bm{\xi}_{t+i})) \right\|^2 \right] \nonumber \\
 &~~~~~+ \EX \left[\left\|\sum_{i=0}^{\ka-1} \bar{\mathbf{Q}} ^{\ka-i-1} \III\left(\n \textbf{F}(\xb_{t+i+1},\bm{\xi}_{t+i+1})-\n \textbf{F}(\xb_{t+i},\bm{\xi}_{t+i})) \right)  \right\|^2\right]  \nonumber \\
&\leq 2\Bigg\{ \EX \left[\left\|\sum_{i=0}^{\ka-1} (\ka-i-1) \bar{\mathbf{Q}} ^{\ka-i-1}\III (\g_{t+i+1}-\g_{t+i}) \right\|^2 \right]+ \EX \left[\left\|\sum_{i=0}^{\ka-1} (\ka-i-1) \bar{\mathbf{Q}} ^{\ka-i-1}\III (\ddd_{t+i+1}-\ddd_{t+i}) \right\|^2\right] \nonumber \\
  &~~~~~+\EX\left[ \left\|\sum_{i=0}^{\ka-1}  (\ka-i-1) \bar{\mathbf{Q}} ^{\ka-i-2}\III(\g_{t+i+1}-\g_{t+i}) \right\|^2\right] + \EX \left[\left\|\sum_{i=0}^{\ka-1} (\ka-i-1) \bar{\mathbf{Q}} ^{\ka-i-2}\III (\ddd_{t+i+1}-\ddd_{t+i}) \right\|^2 \right] \nonumber \\ 
  &~~~~~+ \EX \left[\left\|\sum_{i=0}^{\ka-1} \bar{\mathbf{Q}} ^{\ka-i-1}\III (\g_{t+i+1}-\g_{t+i}) \right\|^2\right]+ \EX \left[\left\|\sum_{i=0}^{\ka-1}  \bar{\mathbf{Q}} ^{\ka-i-1}\III (\ddd_{t+i+1}-\ddd_{t+i}) \right\|^2 \right]\Bigg\}, 
\end{align}
where the inequality is obtained by adding and subtracting
the terms $\n \textbf{f}(\xb_{t+i+ 1}),\,\n\textbf{f}(\xb_{t+i})$, and $\n \textbf{f}(\xb^*)$ in each of the three terms in the first equality. We bound the first term on the right hand side of \eqref{sixterms} and follow a similar approach to bound the rest of the terms. However, before proceeding, we state the following fact whose proof is provided at the end of this appendix:
\begin{equation}\label{proof_app}
   \left\|(i+1)\bar{\textbf{Q}}^{i+1}\III -i \bar{\textbf{Q}}^{i}\III\right\|^2\leq 4, \quad \forall i \in \N.
\end{equation}
The first term on the right hand side of \eqref{sixterms} is bounded as,
\begin{align}\label{st1}
   \EX  &\left[ \left\|\sum_{i=0}^{\ka-1} (\ka- i-1)\bar{\mathbf{Q}}^{\ka- i-1}\III(\g_{t+i+1}-\g_{t+i}) \right\|^2  \right]  \nonumber \\ 
   &~~~~= \EX \Bigg[\Bigg\|  \sum_{i=1}^{\ka-1} \big((\ka - i) \bar{\mathbf{Q}}^{\ka-i}\III - (\ka-i-1)\bar{\mathbf{Q}}^{\ka-i-1}\big) \g_{t+i} -(\ka-1)  \bar{\mathbf{Q}}^{\ka-1}\III \g_{t}\Bigg\|^2 \Bigg]  \nonumber \\
   &~~~~~\leq  
   \sum_{i=1}^{\ka-1} \left\|(\ka - i)  \bar{\mathbf{Q}}^{\ka-1}\III - (\ka-i-1)\bar{\mathbf{Q}}^{\ka-i-1}\III\right\|^2 \EX [\|\g_{t+i}\|^2] 
   + \left\| (\ka-1)\bar{\mathbf{Q}}^{\ka-1}\III \right\|^2 \EX  [\|\g_{t}\|^2]\nonumber \\
    &~~~~~\leq  4\sum_{i=1}^{\ka-1} \EX [\| \g_{t+i}\|^2] + n\sigma_g^2  \leq  4\ka n\sigma_g^2  
\end{align}
where the first inequality is due to Assumption~\ref{asmp4} and the fact that the cross terms of the form $\EX[\langle \g_{p}, \g_{p'} \rangle] = 0$, for any $p<p'$, and the second the inequality is due to Assumption~\ref{asmp4}, \eqref{proof_app}, and the fact that $\ka \| \II-\gamma \Q\|^{\ka-1}\|\III\|\leq  1$(cf. (\ref{Bmatrix})). 
Following a similar approach, we can bound the rest of the terms involving $\textbf{g}_k$ as:
\begin{align}\label{st1-t}
   \EX  \left[ \left\|\sum_{i=0}^{\ka-1} (\ka- i-1)\bar{\mathbf{Q}}^{\ka- i-2}\III(\g_{t+i+1}-\g_{t+i}) \right\|^2  \right]  &\leq 4\ka n\sigma_g^2 \nonumber \\
   \EX  \left[ \left\|\sum_{i=0}^{\ka-1}\bar{\mathbf{Q}}^{\ka- i-1}\III(\g_{t+i+1}-\g_{t+i}) \right\|^2  \right]  &\leq  4(\ka+1) n\sigma_g^2.  
\end{align}
Similarly, considering the second term in (\ref{sixterms}), we have,
\begin{align}\label{intd_t}
   \EX &\left[\left\|\sum_{i=0}^{\ka-1} (\ka- i-1)\bar{\mathbf{Q}}^{\ka- i-1}\III(\ddd_{t+i+1}-\ddd_{t+i}) \right\|^2\right] \nonumber \\
   &~~~~~\leq\ka\left( (\ka-1)  \|\bar{\mathbf{Q}}^{\ka-1}\III\|^2 \EX [\left\|\ddd_{t}\|^2\right] 
    + \sum_{i=1}^{\ka-1} \left\|(\ka - i) \bar{\mathbf{Q}}^{\ka-1}\III - (\ka-i-1)\bar{\mathbf{Q}}^{\ka-i-1}\III\right\|^2 \EX [\|\ddd_{t+i}\|^2]\right) \nonumber\\
    &~~~~~\leq 4\ka\sum_{i=0}^{\ka-1}\EX [ \left\| \ddd_{t+i}\right\|^2 ],
\end{align}
where the first inequality is due to the fact that $\left\|\sum_{i=0}^{\tau-1}a_{i}\right\|^2 \leq \tau \sum_{i=0}^{\tau-1} \|a_i\|^2$ for any $a \in \mathbb{R}^d$.
The same bound also holds for the fourth term in (\ref{sixterms}) while for the last term, we have,
\begin{align}\label{intd}
   \EX \left[\left\|\sum_{i=0}^{\ka-1}\bar{\mathbf{Q}}^{\ka- i-1}\III(\ddd_{t+i+1}-\ddd_{t+i}) \right\|^2\right]  &\leq 4\ka\sum_{i=0}^{\ka}\EX [ \left\| \ddd_{t+i}\right\|^2 ].
\end{align}
We next bound the summation $\sum_{i=0}^{\ka}\EX [ \left\| \ddd_{t+i}\right\|^2$. For all $i<\tau$:
\begin{align}\label{ll012}
\|\ddd_{t+i}\|^2 &= \| \n \textbf{f}(\xb_{t+i}) - \n \textbf{f}(\bar{\xb}_{t+i})+\n \textbf{f}(\bar{\xb}_{t+i})-\n \textbf{f}(\textbf{x}^*) \|^2\nonumber\\
& \leq 2L^2 \|\Psi_{t+i} \|^2 +2\|\n \textbf{f}(\bar{\textbf{x}}_{t+i})-\n \textbf{f}(\textbf{x}^*) \|^2,
\end{align}
where the second inequality is due to Assumption \ref{asmp3}. Now, for $i=\tau$, we have,  
\begin{align}\label{ttaubd}
  \|\ddd_{t+\tau}\|^2 &\leq  2L^2\|\Psi_{t+\ta} \|^2 +2\|\n \textbf{f}(\bar{\textbf{x}}_{t+\ka})-\n \textbf{f}(\textbf{x}^*) \|^2\nonumber\\
   &\leq  2L^2\|\Psi_{t+\ta} \|^2+4\|\n \textbf{f}(\bar{\textbf{x}}_{t+\ka})-\n \textbf{f}(\bar{\textbf{x}}_{t+\ka-1}) \|^2 +4 \|\n \textbf{f}(\bar{\textbf{x}}_{t+\ka-1})-\n \textbf{f}(\textbf{x}^*) \|^2 \nonumber\\
   &\leq  2L^2\|\Psi_{t+\ta} \|^2+4L^2\|
 \bar{\textbf{x}}_{t+\ka}- \bar{\textbf{x}}_{t+\ka-1} \|^2 +4 \|\n \textbf{f}(\bar{\textbf{x}}_{t+\ka-1})-\n \textbf{f}(\textbf{x}^*) \|^2. 
\end{align}
The expression for $\bar{\textbf{x}}_{t+\ka}$ can be written as (cf. \eqref{ll}),
\begin{align*}
\xx_{t+\ka} &= \bar{\textbf{x}}_{t+\ka-1}-   \frac{\alpha}{n} \big(1_n1_n^T\otimes I_d \big) \n \textbf{f} (\xb_{t+\tau-1})  + \alpha \bm{\epsilon}_{t+\tau-1,g} + \gamma \bar{\bm{\epsilon}}_{t+\tau-1,c} - \alpha\gamma \sum_{j=0}^{t+\tau-2}\bar{\bm{\epsilon}}_{j,c},
\end{align*}
where $\bm{\epsilon}_{k,g} \defeq \frac{1}{n} \big(1_n1_n^T\otimes I_d \big)\Big( \n \textbf{f}( \xb_{k}) -  \n \textbf{F} (\xb_{k},\bm{\xi}_k)\Big)$ and $\bar{\bm{\epsilon}}_{k,c} \defeq \frac{1}{n} \big(1_n1_n^T\otimes I_d \big)\hat{\textbf{Q}} \bm{\epsilon}_{k,c}$ for any $k \in \N$. Taking square norms and expectations, we get,
\begin{align}\label{ll1}
\EX &[\|\bar{\textbf{x}}_{t+\tau} - \bar{\textbf{x}}_{t+\tau-1} \|^2 | \mathcal{F}_{t+\tau-1}] \nonumber \\
&~~~~~~\leq  2\EX \left[\left\|  \frac{\alpha}{n} \big(1_n1_n^T\otimes I_d \big) \left(\n \textbf{f} (\xb_{t+\tau-1})  - \n \textbf{f}(\xb^*) \right)  + \alpha \bm{\epsilon}_{t+\tau-1,g} + \gamma \bar{\bm{\epsilon}}_{t+\tau-1,c}\right\|^2 \Big| \mathcal{F}_{t+\tau-1}\right]\nonumber \\
&\qquad \qquad +  2\alpha^2\gamma^2 \EX\left[\left\| \sum_{j=0}^{t+\tau-2}\bar{\bm{\epsilon}}_{j,c}\right\|^2\Big|  \mathcal{F}_{t+\tau-1}\right],
\end{align}
where we used the fact that $\frac{1}{n}\left(1_n1_n^T\otimes I_d \right) \n \textbf{f}(\textbf{x}^*) = 0$. From Assumptions \ref{asmp2} and \ref{asmp4}, we have for all $k\in \mathbb{N}$,
\begin{equation}\label{del1-t}
 \EX[ \alpha \bm{\epsilon}_{k,g} + \gamma \bar{\bm{\epsilon}}_{k,c} | \mathcal{F}_k] =0,  \qquad \qquad   \EX[\left\|\alpha \bm{\epsilon}_{k,g} + \gamma \bar{\bm{\epsilon}}_{k,c}\right\|^2]\leq\alpha^2 \sigma_g^2 + \gamma^2\sigma^2_{c}, 
\end{equation}
and
\begin{align}\label{epsbound-t}
  \EX \left[\left\| \sum_{j=0}^{k-1}\bar{\bm{\epsilon}}_{j,c} \right\|^2\right] &=  \EX \left[ \sum_{j=0}^{k-1}\left\|\bar{\bm{\epsilon}}_{j,c} \right\|^2\right]  + \sum_{1\leq p,p' \leq k-1}\EX \left[\langle \bar{\bm{\epsilon}}_{p,c} ,\bar{\bm{\epsilon}}_{p'
  ,c} \rangle \right]  
  \leq \sum_{j=0}^{k-1} \sigma_c^2 =  k\sigma_c^2,
\end{align}
where the last inequality is due to the fact that
$ \EX [\langle \bar{\bm{\epsilon}}_{p,c}, \bar{\bm{\epsilon}}_{p',c}\rangle] =  \EX [ \EX  \left[\langle \bar{\bm{\epsilon}}_{p,c}, \bar{\bm{\epsilon}}_{p',c}\rangle|\mathcal{F}_{p'}\right]]=0$ for any $p<p'$. Combining \eqref{ll1}, \eqref{del1-t} and \eqref{epsbound-t}, we have,
\begin{align}\label{ll1-t}
\EX &[\|\bar{\textbf{x}}_{t+\tau} - \bar{\textbf{x}}_{t+\tau-1} \|^2]\nonumber \\
&~~~~~\leq  2\EX\left[\left\|  \frac{\alpha}{n} \big(1_n1_n^T\otimes I_d \big) \left(\n \textbf{f} (\xb_{t+\tau-1})  - \n \textbf{f}(\xb^*) \right)\right\|^2 \right] +2\alpha^2 \sigma_g^2 +  2\left( 1 +  \alpha^2 (t+\tau)\right)\gamma^2 \sigma^2_{c} \nonumber\\
&~~~~~\leq  2\EX\left[\left\|  \frac{\alpha}{n} \big(1_n1_n^T\otimes I_d \big) \left(\n \textbf{f} (\xb_{t+\tau-1}) - \n \textbf{f}(\bar{\textbf{x}}_{t+\ka-1})+\n \textbf{f}(\bar{\textbf{x}}_{t+\ka-1})  - \n \textbf{f}(\xb^*) \right)\right\|^2 \right]\nonumber\\
&\qquad \qquad \qquad \qquad\qquad \qquad\qquad \qquad \qquad \qquad\qquad \qquad
+2\alpha^2 \sigma_g^2 +  2\left( 1 +  \alpha^2 (t+\tau)\right)\gamma^2 \sigma^2_{c} \nonumber\\
&~~~~~\leq  4 \alpha^2 L^2\EX \left[\|\Psi_{t+\tau-1}\|^2\right] + 4\alpha^2\EX\left[ \|\n \textbf{f}(\bar{\textbf{x}}_{t+\ka-1})-\n \textbf{f}(\textbf{x}^*) \|^2\right]   + 2\alpha^2 \sigma_g^2 +  2\left( 1 +  \alpha^2 (t+\tau)\right)\gamma^2 \sigma^2_{c}. 
\end{align}
Taking expectations in (\ref{ttaubd}) and using (\ref{ll1-t}), we get, 
\begin{align}\label{ll123}
\EX  [\|\ddd_{t+\tau}\|^2]  &\leq  2L^2\EX [\|\Psi_{t+\ta} \|^2 ]+4 \EX[\|\n \textbf{f}(\bar{\textbf{x}}_{t+\ka-1})-\n \textbf{f}(\textbf{x}^*) \|^2] +  16L^4\alpha^2\EX[\left\|\Psi_{t+\ta-1} \right\|]^2\nonumber\\ 
& \qquad  + 16\alpha^2 L^2\EX [\left\| \n \textbf{f} (\bar{\xb}_{t+\tau-1})  - \n \textbf{f}(\xb^*) \right\|^2] + 8L^2\left(\alpha^2\sigma^2_{g} + \gamma^2 (1+\alpha^2 (t+\tau)) \sigma_c^2 \right)\nonumber\\
 &\leq 2L^2\EX[ \|\Psi_{t+\ta} \|^2   ] + 5 \EX[\|\n \textbf{f}(\bar{\textbf{x}}_{t+\ka-1})-\n \textbf{f}(\textbf{x}^*) \|^2] + L^2 \EX[\left\|\Psi_{t+\ta-1} \right\|^2 ] \nonumber\\
 &\qquad + 8L^2\left(\alpha^2\sigma^2_{g} + \gamma^2 (1+\alpha^2 (t+\tau)) \sigma_c^2 \right), 
\end{align}
where the last inequality is due to $\alpha^2<1/16L^2$.
Using (\ref{ll012}) and (\ref{ll123})  in (\ref{intd}), we have,
\begin{align}\label{st2}
   \EX \left[\left\|\sum_{i=0}^{\ka-1}\bar{\mathbf{Q}}^{\ka- i-1}\III(\ddd_{t+i+1}-\ddd_{t+i}) \right\|^2\right]  &\leq 12 \ka L^2\sum_{i=0}^{\ka-1}\EX [\| \Psi_{t+i}\|^2] +28\ka\sum_{i=0}^{\ka-1}\EX [\left\|\n \textbf{f}(\bar{\textbf{x}}_{t+i})-\n \textbf{f}(\textbf{x}^*) \right\|^2] \nonumber\\
&\qquad +8\tau L^2\EX [\|\Psi_{t+\ta} \|^2]  + 32\tau L^2\left(\alpha^2\sigma^2_{g} + \gamma^2 (1+\alpha^2 (t+\tau)) \sigma_c^2 \right)
\end{align}
The rest of the terms involving $d_k$ in \eqref{sixterms} can be bounded in the same manner. 
Using \eqref{st1}, \eqref{st1-t} and \eqref{st2} in \eqref{sixterms}, we get,
\begin{align}\label{e2}
 \sum_{i=0}^{\ka-1}\EX \left[ \|\textbf{J}^{\ka-i-1}_\gamma  \textbf{E}^g_{t+i}\|^2\right] &\leq 24n(\ka+1)\sigma^2_{g}  + 72 \ka L^2  \sum_{i=0}^{\ka-1} \EX [\|\Psi_{t+i}\|^2] + 168 \ka \sum_{i=0}^{\ka-1} \EX [\left\|\n \textbf{f}(\bar{\textbf{x}}_{t+i})-\n \textbf{f}(\textbf{x}^*) \right\|^2]\nonumber\\
&\qquad + 48\tau L^2\EX [\|\Psi_{t+\ta} \|^2]  + 192\tau L^2\left(\alpha^2\sigma^2_{g} + \gamma^2 (1+\alpha^2 (t+\tau)) \sigma_c^2 \right).
 \end{align}
Using (\ref{bd1}) and (\ref{e2}) to bound the right hand side in (\ref{decomp1}), we have,
\begin{align}\label{finalbd_t}
&\EX \left[\Big\| \sum_{i=0}^{\ka-1} \textbf{J}^{\ka-i-1}_\gamma \textbf{E}_{t+i} \Big\|^2 \right] \nonumber \\
&~~~~~\leq \frac{16\gamma^2}{\alpha^2} \big(1 + \alpha^2 (\ka^2 +1/2)\big)n\sigma_{c}^2 \ka+48n(\ka+1)\sigma^2_{g}  + 144 \ka L^2  \sum_{i=0}^{\ka-1} \EX [\|\Psi_{t+i}\|^2] \nonumber \\
&~~~~~\qquad + 336 \ka \sum_{i=0}^{\ka-1} \EX \left[\left\|\n \textbf{f}(\bar{\textbf{x}}_{t+i})-\n \textbf{f}(\textbf{x}^*) \right\|^2\right]+
  96\tau L^2\EX [\|\Psi_{t+\ta} \|^2]  + 384\tau L^2\left(\alpha^2\sigma^2_{g} + \gamma^2 (1+\alpha^2 (t+\tau)) \sigma_c^2 \right) \nonumber \\
  &~~~~~\leq 144 \ka L^2  \sum_{i=0}^{\ka-1} \EX [\|\Psi_{t+i}\|^2] 
+ 336 \ka \sum_{i=0}^{\ka-1} \EX \left[\left\|\n \textbf{f}(\bar{\textbf{x}}_{t+i})-\n \textbf{f}(\textbf{x}^*) \right\|^2\right]+
  \frac{1}{4\alpha^2}\EX [\|\Psi_{t+\ta} \|^2]   \nonumber \\
  &~~~~~\qquad + \frac{16\gamma^2}{\alpha^2} \big(2 + \alpha^2 (\ka^2 +1/2) + \alpha^2(t+\tau)\big)n\sigma_{c}^2 \ka+49n(\ka+1)\sigma^2_{g}
\end{align}
where we used $\alpha^2 < 1/384\tau L^2$. Next, taking square norms and expectations in \eqref{e0}, we get, 
\begin{align}\label{finalbd0_t}
\EX \left[ \|\Psi_{t+\ka}\|^2 \right]&= \EX \left[\left\| \textbf{J}^\ka_\gamma \Psi_t + \alpha \sum_{i=0}^{\ka-1} \textbf{J}^{\ka-i-1}_\gamma  \textbf{E}_{t+i}\right\|^2 \right]\nonumber\\
&\leq 2 \EX \left[\left\| \textbf{J}^\ka_\gamma \Psi_t\right\|^2\right] + 2\alpha^2 \EX \left[\left\| \sum_{i=0}^{\ka-1} \textbf{J}^{\ka-i-1}_\gamma  \textbf{E}_{t+i}\right\|^2 \right].
\end{align}
From Lemma~\ref{lem3}, it follows that there there exists a $\ka$ such that $\|\textbf{J}^\ka_\gamma\|^2 \leq \frac{1}{4}\rho'$ for a given $\rho'\in(0,1/4]$. Therefore, 
\begin{equation}\label{rho'eq}
4\| \textbf{J}^\ka_\gamma \Psi_t\|^2 \leq 4\|\textbf{J}^\ka_\gamma\|^2 \| \Psi_{t} \|^2  \leq \rho'  \| \Psi_{t} \|^2 .
\end{equation}
To conclude, we substitute \eqref{finalbd_t} and \eqref{rho'eq} in \eqref{finalbd0_t} to get the required inequality,
\begin{align*}
\EX \left[\|\Psi_{t+\ka}\|^2 \right] &\leq  2\EX\left[\left\| \textbf{J}^\ka_\gamma \Psi_t\right\|^2\right]+ 288 \alpha^2 \ka L^2  \sum_{i=0}^{\ka-1} \EX \left[\|\Psi_{t+i}\|^2\right] 
+ 672  \alpha^2 \ka \sum_{i=0}^{\ka-1} \EX \left[\left\|\n \textbf{f}(\bar{\textbf{x}}_{t+i})-\n \textbf{f}(\textbf{x}^*) \right\|^2\right]+
 \frac{1}{2}\EX [\|\Psi_{t+\ta} \|^2]   \\
 &\qquad + 32\gamma^2 \big(2 + \alpha^2 (\ka^2 +1/2) + \alpha^2(t+\tau)\big)n\sigma_{c}^2 \ka+98n(\ka+1)\alpha^2\sigma^2_{g}\\
 \EX \left[\|\Psi_{t+\ka}\|^2 \right]  &\leq \rho' \EX \|\Psi_{t}\|^2+ 576 \alpha^2 \ka L^2  \sum_{i=0}^{\ka-1} \EX [\|\Psi_{t+i}\|^2] 
+ 1344  \alpha^2 \ka \sum_{i=0}^{\ka-1} \EX \left[\left\|\n \textbf{f}(\bar{\textbf{x}}_{t+i})-\n \textbf{f}(\textbf{x}^*) \right\|^2\right]  \\
  &\qquad + 64\gamma^2 \big(2 + \alpha^2 (\ka^2 +1/2) + \alpha^2(t+\tau)\big)n\sigma_{c}^2 \ka+196n(\ka+1)\alpha^2\sigma^2_{g},
\end{align*}
which proves the bound \eqref{main_rec}. 
The bound (\ref{main_rec_2}) for $\ell<\tau$ is proved exactly along the same lines with the only modification being that the first term is scaled by $\|\textbf{J}^\ell_\gamma\|^2,\,\ell<\tau$ instead of $\|\textbf{J}^\tau_\gamma\|^2$. The former can be bounded by using the expression for $\textbf{J}_\gamma^\ell$ (cf. \ref{matproddef}) as follows:
\begin{align*}
    \|\textbf{J}_\gamma^\ell \Psi_0\|^2 &\leq \big\|(\III-\gamma \Q)^\ell -\ell (\III-\gamma \Q)^\ell\big\|^2\|\Delta \textbf{v}_0\|^2 + \big\|(\III-\gamma \Q)^\ell -\ell (\III-\gamma \Q)^{\ell-1}\big\|^2\|\Delta \textbf{x}_0\|^2 \nonumber \\
    &\qquad + \alpha^2 \left\|\left(\III-\gamma \Q\right)^\ell\right\|^{2}\|\Delta \textbf{y}_0\|^2 \nonumber\\
    &\leq  2(1+\ell^2) (\|\Delta \textbf{v}_0\|^2 + \|\Delta \textbf{x}_0\|^2 + \alpha^2\|\Delta \textbf{y}_0\|^2) \nonumber\\
    &\leq   2(1+\ka^2)\|\Psi_0\|^2
\end{align*}
where the second inequality is due to $ \left\|\left(\III-\gamma \Q\right)^\ell\right\|^{2}\leq 1$ and the last inequality is due to $\ell < \ka$.
\end{proof}
\noindent To conclude, we provide the proof of (\ref{proof_app}).\\

\noindent \textit{Claim:} $\|(i+1)\bar{\mathbf{Q}}^{i+1}\III -i \bar{\mathbf{Q}}^{i}\III\|^2\leq4$.
\begin{proof}
We have
\begin{align*}
    \|(i+1)\bar{\mathbf{Q}}^{i+1}\III -i \bar{\textbf{Q}}^{i}\III\|^2 &\leq \|(i+1)\bar{\mathbf{Q}}^{i+1} -i \bar{\textbf{Q}}^{i}\|^2\|\III\|^2 \\
    &\leq \max_{j\in [n]} |(i+1)(1-\gamma(1-\lambda_j)^{i+1} -  i(1-\gamma(1-\lambda_j)^{i} |^2\\
    &= \big|(1+i)(1-\gamma(1-\la))^{i+1} -  i(1-\gamma(1-\la))^{i} \big|^2 \,\,(\text{for some } \la)\\
    &= \big|(1-\gamma(1-\la))^{i+1}  + i (1-\gamma(1-\la))^{i} (1- \gamma(1-\la)-1)\big|^2\\ 
     &= |(1-\gamma(1-\la))^{i+1}  -i \gamma(1-\la) (1-\gamma(1-\la))^{i} |^2\\
     &\leq 2|(1-\gamma(1-\la))|^{2i+2}  +2 \gamma^2(1-\la)^2 \big(\underbrace{i(1-\gamma(1-\la))^{i}}_{\leq \frac{1}{\gamma(1-\la)}}\big)^2 \\
    &\leq 4
\end{align*}
where the second inequality is due to $\|\III\| \leq 1$ and the last inequality is due to $\gamma(1-\la) \in [0,1]$.
\end{proof}

\section*{Appendix IV: Proof of Lemma \ref{lem4}}\label{sec:apndxIV}
\begin{proof}
We begin by defining the following quantities, for any $t\geq\ta$,
\begin{equation}\label{EP}
    A_{t} \defeq\frac{1}{\ka} \sum_{i=t-\ta}^{t-1} a_{i} \quad \text{and} \quad E_{t} \defeq \sum_{i=t-\ta}^{t-1} e_{i}.
\end{equation}
For future reference, we note that for the index $t=k+j$ with $j<\ta \leq t$, \eqref{rel} can be expressed in terms of $A_{k+j}$ and $E_{k+j}$
\begin{equation}\label{alt_ineq}
a_{k+j} \leq \rho' a_{k+j-\ka} + bA_{k+j} + c   E_{k+j} +r
\end{equation}
\textbf{Step (i):} We first prove a recursive relation for $A_{k+\ta}$ in terms of $A_k$ for $k\geq \ta$  and  $E_{k+i}$ for $0\leq i\leq \ka-1$. We begin by considering $A_{k+j}$ for any $j<\ta$,
\begin{align*}
A_{k+j} &= \frac{1}{\ka}\left((a_{k+j-\ta} +\cdots a_{k-1})+ (a_{k} +\cdots + a_{k+j-1})\right)  \nonumber \\
&= \frac{1}{\ta} \left(\sum_{i=j}^{\ka-1} a_{k+i-\ta}  + \sum_{i=0}^{j-1} a_{k+i}\right).
\end{align*}
By \eqref{alt_ineq} (with $j=i$), substituting for $a_{k+i}$ in the second summation above, 
\begin{align*}
    A_{k+j}&\leq \frac{1}{\ka}\left( \sum_{i=j}^{\ka-1} a_{k+i-\ta}   + \sum_{i=0}^{j-1} \rho' a_{k+i-\ka}\right) +\frac{b}{\ka}\sum_{i=0}^{j-1} A_{k+i} +\frac{c}{\ka} \sum_{i=0}^{j-1} E_{k+i} +r \\
    &\leq  \frac{1}{\ka}\left(a_{k-1} +\cdots+ a_{k+j-\ta} + a_{k+j-\ta-1} + a_{k-\ta} \right)+ \frac{b}{\ka}\sum_{i=0}^{j-1} A_{k+i} + \frac{c}{\ka} \sum_{i=0}^{j-1} E_{k+i} +r \\
    &=  A_k + \frac{b}{\ka}\sum_{i=0}^{j-1} A_{k+i} + \frac{c}{\ka} \sum_{i=0}^{j-1} E_{k+i} +r,
\end{align*}
where the second inequality holds since $\rho'<1$. Thus, if follows for $j<\ta$,
\begin{equation}\label{tt1}
    A_{k+j} \leq  A_{k} + \frac{b}{\ka} \sum_{i=0}^{j-1}A_{k+i} + \frac{c}{\ka}   \sum_{i=0}^{j-1} E_{k+i} +r.
\end{equation}
By the definition of $A_{k+\ta}$ \eqref{EP} and \eqref{alt_ineq},
\begin{equation}\label{tt}
\begin{aligned}
    A_{k+\ka} =\frac{1}{\ka}\sum_{j=0}^{\ka-1} a_{k+j}&\leq \frac{1}{\ka}\sum_{j=0}^{\ka-1} \left(\rho'a_{k+j-\ka}  + b  A_{k+j} + c E_{k+j}   +r \right)\\
    &= \rho' A_k + \frac{b}{\ka} \sum_{j=0}^{\ka-1} A_{k+j} + \frac{c}{\ka} \sum_{j=0}^{\ka-1} E_{k+j}  +r\\
    &= \rho' A_k +\frac{b}{\ka} A_{k+\ta-1}+ \frac{b}{\ka} \sum_{j=0}^{\ka-2} A_{k+j} + \frac{c}{\ka} \sum_{j=0}^{\ka-2} E_{k+j} +\frac{c}{\ka} E_{k+\ta-1} +r,  
\end{aligned}
\end{equation}
Next, by \eqref{tt1} with $j=\ta-1$ and $b < \rho'/4$, it follows that
\begin{align*}
    A_{k+\ka} &\leq  \rho' A_k +\frac{b}{\ka} A_{k+\ta-1}+ \frac{b}{\ka} \sum_{j=0}^{\ka-2} A_{k+j} + \frac{c}{\ka} \sum_{j=0}^{\ka-2} E_{k+j} +\frac{c}{\ka} E_{k+\ta-1} + r\\
    &\leq \rho' A_k
  + \frac{b}{\ka}\left(A_{k} +\frac{b}{\ka} \sum_{j=0}^{\ka-2}A_{k+j} + \frac{c}{\ka}  \sum_{j=0}^{\ka-2} E_{k+j} + r\right)\\
    &   \qquad \qquad+ \frac{b}{\ka}\sum_{j=0}^{\ka-2} A_{k+j} + \frac{c}{\ka} \sum_{j=0}^{\ka-2} E_{k+j}   + \frac{c}{\ka}E_{k+\ka-1} + r \\
  &\leq   \rho'\left(1+\frac{1}{4\ka}\right) A_k
 + \left(1+\frac{b}{\ka}\right) \left( \frac{b}{\ka}\sum_{j=0}^{\ka-2}A_{k+j} + \frac{c}{\ka} \sum_{j=0}^{\ka-2} E_{k+j}  + r \right) +\frac{c}{\ka} E_{k+\ka-1}.
\end{align*}
Recursive application of the above, over $j$ for $1\leq j \leq \ka -1$ yields the following inequality
\begin{equation}\label{majineq}
\begin{aligned}
  A_{k+\ka} &\leq   \rho'\Bigg(1+\frac{1}{4\ka}\Bigg)^{\ka-1} A_k
 + \frac{c}{\ka}  \sum_{j=0}^{\ka-1} \Big(1+\frac{b}{\ka}\Big)^{\ka-j-1} E_{k+j} +  \Big(1+\frac{b}{\ka}\Big)^{\ka-1} r \\
  &\leq   (1-\rho) A_k
 + \frac{2c}{\ka} \sum_{j=0}^{\ka-1}  E_{k+j} +2r,   
 \end{aligned}
 \end{equation}
where $b\leq 1/4$, 
$\left(1+\frac{b}{\ka}\right)^p \leq \left(1+\frac{1}{4\ka}\right)^p   \leq \exp(1/4) \leq 2
$ for any $p\leq \ka$ and  $2\rho'= 1-\rho $. Next, recall \eqref{alt_ineq} with $j=\ta$:
\begin{equation}\label{eq.1}
\begin{aligned}
    a_{k+\ka} &\leq  \rho' a_k + b A_{k+\ka} + c  \sum_{j=0}^{\ka-1} e_{k+j} + r \\
  & \stackrel{(\ref{majineq})}{ \leq}   \rho' a_k + b(1-\rho) A_k +   \frac{2bc}{\ka} \sum_{j=0}^{\ka-1}  E_{k+j}  + c \sum_{j=0}^{\ka-1}  e_{k+j} +(2b+1)r \\
    &\leq \frac{1}{2}(1-\rho)a_k + \frac{\rho}{4} (1-\rho)A_k + c\sum_{j=0}^{\ka-1}  \Big(\frac{1}{\ka} E_{k+j} + e_{k+j}\Big) +2r, 
\end{aligned}
\end{equation}
where we have used the fact that $\rho' = (1-\rho)/2$, $b\leq\rho'/4<\rho/4$ for $\rho'\in (0,1/4]$ to get the last inequality. Adding, \eqref{majineq} and \eqref{eq.1}, 
\begin{equation}\label{tmp3}
\begin{aligned}
    A_{k+\ka} + a_{k+\ka} &\leq (1-\rho)\left(1+\frac{\rho}{4}\right)(A_k + a_{k}) +  \frac{3c}{\ka}\sum_{j=0}^{\ka-1}  E_{k+j} + c\sum_{j=0}^{\ka-1} e_{k+j} +4r\\
    &\leq       \left(1-\frac{3\rho}{4}\right)(A_k + a_{k}) + c \sum_{j=0}^{\ka-1} \left(\frac{3}{\ka} E_{k+j}+ e_{k+j}\right) +4r. 
\end{aligned}
\end{equation}

\textbf{Step (ii):} In this step, we establish a descent relation for $A_{t}+a_{t}$. With $k=(m-1)\ta$ for any integer $m\geq2$, \eqref{tmp3} can be expressed as
\begin{equation}\label{mn1}
\begin{aligned}
    A_{m\ta} + a_{m\ta} &\leq \left(1-\frac{3\rho}{4} \right) \left(A_{(m-1)\ka} + a_{(m-1)\ka} \right) + c \sum_{i=0}^{\ka-1} \left(\frac{3}{\ka} E_{(m-1)\ka+i}+ e_{(m-1)\ka+i}\right)+4r \\
    &\leq \left(1-\frac{3\rho}{4} \right)^{m-1} \left(A_{\ta} + a_{\ta} \right) + c \sum_{j=1}^{m-1}   \left(1-\frac{3\rho}{4}\right)^{(m-j-1)} \sum_{i=0}^{\ka-1} \left(\frac{3}{\ka} E_{j\ka+i}+ e_{j\ka+i}\right)\\
    & + 4r \sum_{j=0}^{m-1}   \left(1-\frac{3\rho}{4}\right)^{j}
\end{aligned}
\end{equation}
Let $t=m\ta$ and $m\geq 2$. To bound the summation in \eqref{mn1}, we note that
\begin{align*}
    \left(1-\frac{3\rho}{4}\right)^{m-j-1} &= \left(1-\frac{3\rho}{4}\right)^{-1}\left(1-\frac{3\rho}{4}\right)^{\frac{t-j\ta}{\ka}}\leq c'\left(1- \frac{3\rho }{4\ka}\right)^{t-j\ta}\label{bb1}
\end{align*}
where $c':=  \left(1-\frac{3\rho}{4}\right)^{-1}$, $t=m\ta$ and $(1-x)^a\leq 1-ax$ for $a,x\in [0,1]$. Thus, it follows that
\begin{equation}\label{tmp10}
\begin{aligned}
    \sum_{j=1}^{m-1}   \left(1-\frac{3\rho}{4}\right)^{(m-j-1)}\sum_{i=0}^{\ka-1} \left(\frac{3}{\ka}E_{j\ka+i}+e_{j\ka+i}\right)  &\leq   c' \sum_{j=1}^{m-1}   \left(1-\frac{3\rho}{4\ta}\right)^{t-j\ta}\sum_{i=0}^{\ka-1}  \left(\frac{3}{\ka}E_{j\ka+i}+e_{j\ka+i}\right)\\
    &\leq c'\sum_{j=1}^{m-1}  \sum_{i=0}^{\ka-1} \Big(1- \frac{3\rho}{4\ka}\Big)^{t-{j\ka-i}}  \left(\frac{3}{\ka}E_{j\ka+i}+e_{j\ka+i}\right)\\
    &\leq  c'  \sum_{p=\ta}^{t-1} \Big(1- \frac{3\rho}{4\ka}\Big)^{t-p}  \left(\frac{3}{\ka}E_{p}+e_p\right),
\end{aligned}
\end{equation}
where $\left(1-3\rho/4\ka\right)^{-i}\geq 1$ since $\rho<1$ and the index  $p:= j\ka+i$. To bound the term involving $E_p$ in the summation in \eqref{tmp10}, it follows by \eqref{EP} 
\begin{align*}
    \Big(1- \frac{3\rho}{4\ka}\Big)^{t-p}  E_{p}\leq  \Big(1- \frac{3\rho}{4\ka}\Big)^{t-p}\sum_{i=p-\ka}^{p-1}e_i &\leq  \Big(1- \frac{3\rho}{4\ka}\Big)^{-\ta}  \sum_{i=p-\ka}^{p-1}\Big(1- \frac{3\rho}{4\ka}\Big)^{t-i}e_i\\
     &\leq  \Big(1+\frac{3\rho}{4\ka}\Big)^{\ka}  \sum_{i=p-\ka}^{p-1}\Big(1- \frac{3\rho}{4\ka}\Big)^{t-i}e_i  \\
     &\leq  3 \sum_{i=p-\ka}^{p-1}\Big(1- \frac{3\rho}{4\ka}\Big)^{t-i}e_i,
\end{align*}
where the second inequality follows due to $(1-x)^{-1}<1+x$ for $x\in (0,1)$ and the last inequality follows due to  $\left(1+3\rho/4\ka\right)^{\ka} \leq \exp(3\rho/4)\leq 3$ for $\rho<1$. Summing the above for $p=\ta$ to $t-1$
\begin{align}\label{e-ineq}
   \sum_{p=\ta}^{t-1} \Big(1- \frac{3\rho}{4\ka}\Big)^{t-p}  E_{p}
    &\leq 3\sum_{p=\ka}^{t-1}  \sum_{i=p-\ka}^{p-1}\Big(1- \frac{3\rho}{4\ka}\Big)^{t-i}e_i\nonumber\\
    &=3\left( \sum_{i=t-\ka-1}^{t-2}\Big(1- \frac{3\rho}{4\ka}\Big)^{t-i}e_i+\sum_{i=t-\ka-2}^{t-3}\Big(1- \frac{3\rho}{4\ka}\Big)^{t-i}e_i  +\cdots+ \sum_{i=0}^{\ta-1}\Big(1- \frac{3\rho}{4\ka}\Big)^{t-i}e_i \right)\nonumber\\
     &\leq  3\ka \sum_{p=0}^{t-2} \Big(1- \frac{3\rho}{4\ka}\Big)^{t-p}e_p.  
\end{align}
Substituting \eqref{e-ineq} into \eqref{tmp10}, 
\begin{equation}\label{eq2_albert}
\begin{aligned}
    \sum_{j=1}^{m-1}   \left(1-\frac{3\rho}{4}\right)^{(m-j-1)}\sum_{i=0}^{\ka-1}  \left(\frac{3}{\ka}E_{j\ka+i}+e_{j\ka+i}\right) 
    &\leq  10c' \sum_{p=0}^{t-1} \left(1- \frac{3\rho}{4\ka}\right)^{t-p}e_p\\
    &\leq  40 \sum_{p=0}^{t-1} \left(1- \frac{3\rho}{4\ka}\right)^{t-p}e_p,
\end{aligned}
\end{equation}
where for $\rho<1$ and $c':=\left(1-3\rho/4\right)^{-1}$ it follows that $c'\leq 4$. Substituting \eqref{eq2_albert} into \eqref{mn1}, if $t=m\ka$ and $m\geq 2$, 
\begin{equation}\label{rel2}
\begin{aligned}
    a_t \leq A_{t} + a_{t} &\leq \Big(1- \frac{3\rho}{4}\Big)^{m-1}(A_{\ta} + a_{\ta}) + 40 c\sum_{j=0}^{t-1}   \Big(1- \frac{3\rho}{4\ka}\Big)^{t-j}e_j + 4r \sum_{j=0}^{m-1}   \left(1-\frac{3\rho}{4}\right)^{j}\\
    &\leq \Big(1- \frac{3\rho}{4}\Big)^{\frac{t}{\ta}-1}(A_{\ta}+a_{\ta}) + 40 c\sum_{j=0}^{t-1}   \Big(1- \frac{3\rho}{4\ka}\Big)^{t-j}e_j+ \frac{16r}{3\rho}, \\
     &\leq \Big(1- \frac{3\rho}{4\ta}\Big)^{t-\ta}(A_{\ta}+a_{\ta}) + 40 c\sum_{j=0}^{t-1}   \Big(1- \frac{3\rho}{4\ka}\Big)^{t-j}e_j +\frac{16r}{3\rho}, 
\end{aligned}
\end{equation}
where we have used $(1-x)^\frac{1}{\ta} <1-\frac{x}{\ta}$ for $x \in [0,1]$ and $\ta\in \mathbb{N}$ to get the last inequality. For the case where $t=m\ta+\ell$, $m\geq 2$ and $\ell<\ta$,  the above bound is
\begin{equation}\label{lastbound}
\begin{aligned}
    a_t &\leq \Big(1- \frac{3\rho}{4\ka}\Big)^{t-\ta-\ell} (A_{\ta+\ell} + a_{\ta+\ell}) + 40c\sum_{j=\ell}^{t-1}   \Big(1- \frac{3\rho}{4\ka}\Big)^{t-j}e_j +\frac{16r}{3\rho}\\
    &\leq 5\Big(1- \frac{3\rho}{4\ka}\Big)^{t} (A_{\ta+\ell} + a_{\ta+\ell}) + 40c\sum_{j=\ell}^{t-1}   \Big(1- \frac{3\rho}{4\ka}\Big)^{t-j}e_j+\frac{16r}{3\rho}, 
\end{aligned}
\end{equation}
where the last inequality holds due to $\left(1- \frac{3\rho}{4\ka}\right)^{-\ell-\ta} <\exp\left( 3\rho/2\right)< 5$ for $\ell<\ta$.

\textbf{Step (iii):} In this step, we  use \eqref{rel} to bound the term $A_{\ta+\ell} + a_{\ta+\ell}$ in \eqref{lastbound} where $\ell<\ta$. The argument is similar to the one employed in \textbf{Step (i)}  with appropriate modifications. By \eqref{rel}, for $t=j$ and $j<\ka$, 
\begin{equation}\label{tmp}
    a_{j} \leq \rho'' a_{0} + \frac{b}{\ka} \sum_{i=0}^{j-1} a_{i} + c \sum_{i=0}^{j-1} e_{i} + r.  
\end{equation}
Note, by \eqref{rel} and $\rho'<1$, \eqref{tmp} holds for $\ta\leq j<2\ta$ with a larger $r$,
\begin{equation}\label{tmp2}
\begin{aligned}
    a_{j}           &\leq\rho' a_{j-\ta} + \frac{b}{\ka}\sum_{i=j-\ta}^{j-1} a_{i} + c\sum_{i=j-\ta}^{j-1} e_{i}  + r \\ 
    &\leq \rho'\left(\rho''a_0   + \frac{b}{\ka} \sum_{i=0}^{j-\ta-1} a_{i} + c \sum_{i=0}^{j-\ta-1} e_{i} + r \right) + \frac{b}{\ka}\sum_{i=j-\ta}^{j-1} a_{i} + c\sum_{i=j-\ta}^{j-1} e_{i}  +r \\ 
    &\leq \rho'' a_{0} + \frac{b}{\ka} \sum_{i=0}^{j-1} a_{i} + c \sum_{i=0}^{j-1} e_{j} + 2r.
\end{aligned}
\end{equation}
Recursive application of \eqref{tmp2}, over $j$ for $1\leq j < 2\ka $ yields
\begin{align}\label{rel00}
    a_j&\leq\left(1 +\frac{b}{\ka} \right)^{j-1} \rho''a_0+c \sum_{i=0}^{j-1} \left(1 +\frac{b}{\ka} \right)^{j-i-1} e_i + 2 \left(1 +\frac{b}{\ka} \right)^{j-1}r \nonumber\\
    &\leq 2 \rho''a_0+2c \sum_{i=0}^{j-1}  e_i +4r,\qquad j<2\ta
\end{align}
where the inequality holds due to the fact that $\left(1 +b/\ka \right)^{p} \leq \exp(pb/\ka) \leq \exp(2b)< 2$ for $p<2\ta$. Next, by \eqref{rel00} we bound $A_{\ka+\ell}$ for $\ta+\ell<2\ta$
\begin{equation}\label{1bound}
\begin{aligned}
    A_{\ta+\ell} = \frac{1}{\ta} \sum_{j=\ell}^{\ta+\ell-1} a_j &\leq \frac{1}{\ta} \sum_{j=\ell}^{\ta+\ell-1} \left(2\rho'' a_{0} +   2c \sum_{i=0}^{j-1} e_{i}+4r \right) \\
    &\leq  \left(2\rho'' a_{0} +   2c\sum_{i=0}^{\ta+\ell-1} e_{i}  +4 r\right) \left(\frac{1}{\ta}\sum_{j=\ell}^{\ta+\ell-1} 1 \right)\\
     &=  2\rho'' a_{0} +   2c \sum_{i=0}^{\ta+\ell-1}e_{i} +4r. 
\end{aligned}
\end{equation}
Finally, adding \eqref{rel00} with $j = \ka + \ell$ and \eqref{1bound} gives 
\begin{align}\label{ft}
   A_{\ka+\ell} + a_{\ka+\ell}  &\leq  4\rho'' a_{0} + 4c  \sum_{i=0}^{\ta+\ell-1}  e_{i} +4r. 
\end{align}
Substituting for $A_{\ka+\ell} + a_{\ka+\ell}$ in \eqref{lastbound} using the above inequality completes the proof for any $t\geq 2\tau$. For $t<2\tau$, we have from (\ref{rel00}) with $j=t$ and the fact that $\left(1 - \frac{3\rho}{4\tau}\right)^{-t} \leq 2$ for any $t<2\tau$, 
\begin{align*}
  a_{t}  &\leq  2\rho'' a_{0} + 2c  \sum_{i=0}^{t-1}  e_{i} +4r\\
  &\leq 4\rho''\left(1 - \frac{3\rho}{4\tau}\right)^{t} a_{0} + 4c  \sum_{i=0}^{t-1} \left(1 - \frac{3\rho}{4\tau}\right)^{t-i} e_{i} + 8\left(1 - \frac{3\rho}{4\tau}\right)^{t}r,
\end{align*}
implying that \eqref{lem5_main} also holds for any $t<2\tau$.
\end{proof}

\bibliographystyle{IEEEtran}
\bibliography{references}

\end{document}